\documentclass[12pt]{amsart} 

\usepackage{amsmath,amsthm,amsfonts,latexsym,euscript,amssymb}

\newtheorem{DF}{Definition}[section]
\newtheorem{LM}[DF]{Lemma}
\newtheorem{PROP}[DF]{Proposition}
\newtheorem{THM}[DF]{Theorem}
\newtheorem{COR}[DF]{Corollary}
\newtheorem{RMK}[DF]{Remark}
\newtheorem{RMKS}[DF]{Remarks}
\newtheorem{PROB}[DF]{Problem}

\newcommand{\bgdf}{\begin{DF}}
\newcommand{\nddf}{\end{DF}}
\newcommand{\bglm}{\begin{LM}}
\newcommand{\ndlm}{\end{LM}}
\newcommand{\bgprop}{\begin{PROP}}
\newcommand{\ndprop}{\end{PROP}}
\newcommand{\bgth}{\begin{THM}}
\newcommand{\ndth}{\end{THM}}
\newcommand{\bgthm}{\begin{THM}}
\newcommand{\ndthm}{\end{THM}}
\newcommand{\bgcor}{\begin{COR}}
\newcommand{\ndcor}{\end{COR}}
\newcommand{\bgrmk}{\begin{RMK}}
\newcommand{\ndrmk}{\end{RMK}}
\newcommand{\bgrmks}{\begin{RMKS}}
\newcommand{\ndrmks}{\end{RMKS}}
\newcommand{\bgprob}{\begin{PROB}}
\newcommand{\ndprob}{\end{PROB}}
\newcommand{\bgeq}{\begin{eqnarray}}
\newcommand{\ndeq}{\end{eqnarray}}
\newcommand{\bgeqq}{\begin{eqnarray*}}
\newcommand{\ndeqq}{\end{eqnarray*}}

\newcommand{\pf}{{\em Proof. }} 

\newcommand{\vv}{\vspace{2mm}\\}

\def\Uq{U_q ({\mathfrak g})}

\numberwithin{equation}{section}

\newcommand{\dfref}[1]{Definition~\ref{#1}}
\newcommand{\lemref}[1]{Lemma~\ref{#1}}
\newcommand{\lmref}[1]{Lemma~\ref{#1}}
\newcommand{\propref}[1]{Proposition~\ref{#1}}

\newcommand{\thmref}[1]{Theorem~\ref{#1}}

\newcommand{\probref}[1]{Problem~\ref{#1}}
\newcommand{\secref}[1]{\S\ref{#1}}
\begin{document}




\title[Simple Compact Quantum Groups]
{\bf
Simple Compact Quantum Groups I
}
\author[Shuzhou Wang]
{\bf Shuzhou Wang}
\thanks{Research supported in part by the National Science Foundation
grant DMS-0096136}


\address{Department of Mathematics, University of Georgia,
Athens, GA 30602, USA
\newline \indent
Fax: 706-542-5907,
Tel: 706-542-0884
}
\email{szwang@math.uga.edu}

\begin{abstract}
The notion of simple compact quantum group is introduced.
As non-trivial (noncommutative and noncocommutative) examples,
the following families of compact quantum groups are shown to be simple:
(a) The universal quantum groups $B_u(Q)$ for $Q \in GL(n, {\mathbb C})$
satisfying $Q \bar{Q} = \pm I_n$, $n \geq 2$;
(b) The quantum automorphism groups $A_{aut}(B, \tau)$ of finite dimensional
$C^*$-algebras $B$ endowed with the canonical trace $\tau$
when $\dim(B) \geq 4$, including the quantum permutation groups
$A_{aut}(X_n)$ on $n$ points ($n \geq 4$);
(c) The standard deformations $K_q$ of simple compact Lie groups $K$
and their twists $K_q^u$, as well as Rieffel's deformation $K_J$.
\end{abstract}

\maketitle

\section{Introduction}
\label{Introd}

The theory of quantum groups saw spectacular breakthroughs in the
1980's when on the one hand Drinfeld \cite{Dr1} and Jimbo \cite{Jimb1} discovered the
quantized universal enveloping algebras of
semisimple Lie algebras based on the work of
the Faddeev school on the quantum inverse scattering method,
and on the other hand Woronowicz \cite{Wor4,Wor5,Wor6}
independently discovered quantum deformations
of compact Lie groups and formulated the axioms for compact quantum groups.
Further work of Rosso \cite{Rosso87a,Rosso90a},  Soibelman and Vaksman, Levendorskii
\cite{VS,Soib1,Lev1} showed that ``compact real forms''
$K_q$ of the Drinfeld-Jimbo quantum groups
and their twists $K_q^u$ are examples of compact quantum groups in the
sense of Woronowicz. Most notable of these is the work of Soibelman \cite{Soib1}
based on his earlier joint work with Vaksman \cite{VS},
in which a general Kirillov type orbit theory of
representations of the quantum function algebras of deformed simple compact Lie groups
was developed using the orbits of dressing transformations
(i.e. symplectic leaves) in Poisson Lie group theory
(see also the monograph \cite{KorogodskiSoib1} for more detailed
treatment).

Starting in his Ph.D. thesis \cite{Wang}, the author of the present
article took a different direction from the above by viewing quantum
groups as intrinsic objects and found in a series of papers
(including \cite{W5} in collaboration with Van Daele) several
classes of compact quantum groups that can not be obtained as
deformations of Lie groups. The most important of these are the
universal compact quantum groups of Kac type $A_u(n)$ and their
self-conjugate counterpart $A_o(n)$ \cite{W1}, the more general
universal compact quantum groups $A_u(Q)$ and their self-conjugate
counterpart $B_u(Q)$ \cite{W5,W5'}, where $Q \in GL(n, {\mathbb
C})$, and the quantum automorphism groups $A_{aut}(B, tr)$ of finite
dimensional $C^*$-algebras $B$ endowed with a tracial functional
$tr$, including the quantum permutation groups $A_{aut}(X_n)$ on the
space $X_n$ of $n$ points \cite{W15}. Further studies of these quantum groups
reveal remarkable properties: (1) According to deep work of
Banica \cite{Banica1,Banica2,Banica7},
the representation rings (also called the fusion rings)
of the quantum groups $B_u(Q)$ (when $Q \bar{Q}$ is a scalar)
are all isomorphic to that of $SU(2)$ (see Th\'eor\`eme 1 in \cite{Banica1}),
and the representation rings of $A_{aut}(B, \tau)$
(when $\dim(B) \geq 4$, $\tau$ being the canonical trace on $B$)
are all isomorphic to that of $SO(3)$
(see Theorem 4.1 in \cite{Banica7}),
and the representation ring of
$A_u(Q)$ is almost a free product of two copies of ${\mathbb Z}$
(see Th\'eor\`eme 1 in \cite{Banica2});
(2) The compact quantum groups $A_u(Q)$ admit ergodic actions on
both finite and infinite injective von Neumann factors
\cite{W14};
(3) The special $A_u(Q)$'s for positive $Q$ and $B_u(Q)$'s
for $Q$ satisfying the property $Q \bar{Q} = \pm I_n$
are classified up to isomorphism using {\em respectively} the eigenvalues of
$Q$ (see Theorem 1.1 in \cite{W17}) and polar decomposition of $Q$ and eigenvalues
of $|Q|$ (see Theorem 2.4 in \cite{W17}), and the
general $A_u(Q)$'s and $B_u(Q)$'s for arbitrary $Q$
have explicit decompositions as free products of the
former special ones (see Theorems 3.1 and 3.3 and Corollaries 3.2 and 3.4 in \cite{W17});
(4) Certain quantum symmetry groups in the
theory of subfactors were found by Banica \cite{Banica2000a,Banica2002}
to fit in the theory of compact quantum groups;
(5) The quantum permutation groups $A_{aut}(X_n)$ admit
interesting quantum subgroups that appear in connection with
other areas of mathematics, such as the quantum automorphism groups
of finite graphs and the free wreath products discovered by Bichon
\cite{Bichon1,Bichon2}.
See also \cite{Bichon3} and \cite{Banica2005a}--\cite{Banica-Bichon-Collins2007}
and the references therein for other interesting
results related to the quantum permutation groups.

The purpose of this article is to initiate a study of
simple compact quantum groups. It focuses on the introduction of
a notion of simple compact quantum groups and first examples.
It is shown that the compact quantum groups
mentioned in the last two paragraphs are simple in generic cases.
The paper is organized as follows.

 In \secref{notion-of-normal}, we recall
 the notion of a normal quantum subgroup $N$ of a compact quantum group
$G$ introduced in \cite{Wang,W1}, on which the main notion of a
simple compact quantum group in this paper depends. We
prove several equivalent conditions
for $N$ to be normal, including one that stipulates that
the quantum coset spaces $G/N$ and $N \backslash G$ are identical.
Further applications of these are contained in \cite{normal}.

In \secref{simple} the notion of simple compact quantum groups is introduced.
In the classical setting, the notion of a simple compact Lie group can be defined
in two ways: one using Lie algebra and the other using the group itself.
Though the universal enveloping algebras of simple Lie groups can be deformed into
the quantized universal enveloping algebras \cite{Dr1,Jimb1},
we have no analog of Lie algebras for general quantum groups. Hence
we formulate the notion of a {\bf simple compact quantum group} using
group theoretical language so that our notion reduces precisely to
the notion of a {\em simple compact Lie group} when the quantum group
is a compact Lie group:

\bgdf
\label{define-simple}
A compact matrix quantum group is called {\bf simple} if it is connected
and has no non-trivial connected normal quantum subgroups and
no non-trivial representations of dimension one.
\nddf

Here a compact quantum group $G$ is called {\bf connected} if the coefficients of every
non-trivial irreducible representation of $G$ generate an infinite dimensional $C^*$-algebra.
In the classical situation,
the fact that a simple compact Lie group has no non-trivial representations of dimension
one is a consequence of the deep Weyl dimension formula.
It is not known if the postulate that a simple compact matrix quantum group
has no non-trivial representations of dimension one follows from the
other postulates in the definition, for we do not have a dimension formula
for irreducible representations of a general simple compact quantum group
except the specific examples studied in this paper.

After preparatory work in \secref{notion-of-normal} and \secref{simple},
the main examples of this paper are studied in \secref{B_u-A_{aut}} and \secref{K_q-K_J}.
Recall \cite{Banica7} that {\em the canonical trace $\tau$}
on a finite dimensional $C^*$-algebra $B$ is the restriction
of the unique tracial state on the algebra $L(B)$ of operators on $B$.
In \secref{B_u-A_{aut}}, we prove that $B_u(Q)$ and $A_{aut}(B, \tau)$ are simple:
\bgth
\label{simpleB_u(Q)1}
(see \thmref{simpleB_u(Q)})
Let $Q \in GL(n, {\mathbb C})$ be such
that $Q \bar{Q} = \pm I_n$ and $n \geq 2$.
Then $B_u(Q)$ is a simple compact quantum group.
\ndth
\bgth
\label{simpleA_{aut}1}
(see \thmref{simpleA_{aut}})
Let $B$ be a finite dimensional $C^*$-algebra with $\dim(B) \geq 4$
and $\tau$ its canonical trace.
Then $A_{aut}(B, \tau)$ is a simple compact quantum group.
\ndth
\noindent
The proofs of these two results
rely heavily on the fundamental work of Banica \cite{Banica1,Banica7} on the
structure of fusion rings (i.e. representative rings) of these quantum groups, as well as
the technical results on the correspondence between Hopf $*$-ideals and
Woronowicz $C^*$-ideals and the reconstruction of a normal
quantum group from the identity in the quotient quantum group,
which are developed in \secref{B_u-A_{aut}} and
are of interest in their own right.

It is also shown in  \secref{B_u-A_{aut}} that the closely related
quantum group $A_u(Q)$ is not simple for any $n$ and any
$Q \in GL(n, {\mathbb C})$ (see \propref{non-simpleA_u(Q)}).

The last section \secref{K_q-K_J} is devoted to
the standard deformations $K_q$ of simple compact Lie groups,
their twists $K_q^u$ \cite{Soib1,Lev1,LS1},
and Rieffel's quantum groups $K_J$ \cite{Rf8},
where $q \in {\mathbb R} \backslash \{ 0 \}$,
$u \in \wedge^2 (i \mathfrak{t})$ and $J$ is an appropriate
skew-symmetric transformation on the direct sum
$ \mathfrak{t} \oplus  \mathfrak{t}$ of 
Cartan subalgebra
$\mathfrak{t}$ of the Lie algebra of $K$:
\bgth
\label{simpleK_qK^u_q1}
(see  \thmref{simpleK_q} and \thmref{simpleK^u_q})
Let $K$ be a connected and simply connected simple compact Lie group.
Then both $K_q$ and its twists $K_q^u$ are simple compact quantum groups.
\ndth
\bgth
\label{simpleK_J1}
(see \thmref{simpleK_J})
Let $K$ be a simple compact Lie group with a toral subgroup $T$ of
rank at least two.
Then $K_J$ is a simple compact quantum group.
\ndth
\noindent
The proofs of \thmref{simpleK_qK^u_q1} and
\thmref{simpleK_J1} make use of the work of
Lusztig and Rosso \cite{Lus88a,Rosso88a} on representations of quantized
universal enveloping algebras, the work of
Soibelman and Levendorskii \cite{Soib1,Lev1,LS1} on quantum
function algebras of  $K_q$ and $K_q^u$, and the work of
Rieffel \cite{Rf8} and the author \cite{W4}
on strict deformations of Lie groups and quantum groups,
as well as the technical results in \secref{B_u-A_{aut}} mentioned earlier.

Classification of simple compact quantum groups
and their irreducible representations up to isomorphism
are two of the main goals in the study of compact quantum groups.
Namely, one would like to develop a theory of simple compact quantum groups
that parallels the Killing-Cartan theory
and the Cartan-Weyl theory for simple compact Lie groups.
To accomplish the first goal,
one must first construct all simple compact quantum groups.
Though we have given several infinite classes of examples
of these in this article,
it should be pointed out that the construction of simple compact quantum
groups is far from being complete. In fact it is fair to say that we are only
at the beginning stage for this task at the moment.
One indication of this is that
{\em all the simple compact quantum groups known so far have
commutative representation rings, and these rings are
order isomorphic to the representation rings of compact Lie groups}
(we call such quantum groups {\bf almost classical}).
The universal compact matrix quantum groups $A_u(Q)$ have a
``very'' noncommutative representation ring, being close to the free product of two
copies of the ring of integers, according to the fundamental work of Banica \cite{Banica2},
where $Q \in GL(n, {\mathbb C})$ are positive, $n \geq 2$. However,
$A_u(Q)$ are not simple quantum groups (see \secref{B_u-A_{aut}}).
Because of their universal property, $A_u(Q)$ should play an important role
in the construction and classification of simple compact quantum groups
with non-commutative representation
rings. A natural and profitable approach seems to be to
study quantum automorphism groups of appropriate quantum spaces
and their quantum subgroups, such as those in
\cite{W15,W14,W16} and the papers of Banica and
Bichon and their collaborators \cite{Banica2000a}--\cite{Bichon3}.
In retrospect, both simple Lie groups and finite simple groups are
automorphism groups, a similar approach for the theory of simple quantum groups
should also play a fundamental role.
\vv
{\bf Convention and Notation}.
We assume that all Woronowicz $C^*$-algebras
(also called Woronowicz Hopf $C^*$-algebras)
considered in this paper to be full
unless otherwise explicitly stated,
since morphisms between quantum groups are
meaningful only for full Woronowicz $C^*$-algebras
cf. \cite{W1,W3}.
For a compact quantum group $G$, $A_G$, or $C(G)$, denote the
underlying Woronowicz $C^*$-algebra and ${\mathcal A}_G$ denotes the associated
canonical dense Hopf $*$-algebra of quantum representative functions on $G$.
Sometimes we also call $A_G$ a compact quantum group, referring to $G$.
See \cite{W1,Wor5}
for more on
other unexplained definitions and notations used in this paper.

\section{The Notion of Normal Quantum Subgroups}
\label{notion-of-normal}

Before making the notion of simple quantum groups precise,
we recall the notion of normal quantum subgroups
(of compact quantum groups) introduced in \cite{Wang,W1}
and study their properties further.
Let $(N, \pi)$ be a quantum subgroup of a compact quantum group
$G$ with surjections $\pi: A_G \longrightarrow A_N$ and
$\hat{\pi}: {\mathcal A}_G \longrightarrow {\mathcal A}_N$.
The quantum group $(N, \pi)$ should be more precisely called a closed
quantum subgroup, but we will omit the word {\em closed } in this
paper, since we do not consider non-closed quantum subgroups.
Define
$$A_{G/N} = \{ a \in A_G | (id \otimes \pi) \Delta (a) = a \otimes 1_N \},
\; \; \; $$
$$
A_{N \backslash G} = \{ a \in A_G | (\pi \otimes id) \Delta (a)
= 1_N \otimes a \}, \; \; \; $$
where $\Delta$ is the coproduct on $A_G$, $1_N$ is the unit of the algebra
$A_N$. We omit the subscript $N$ in $1_N$ when no confusion arises.
Similarly, we define
$$
{\mathcal A}_{G/N} = {\mathcal A}_G \cap
A_{G/N}, \; \; \;
\text{and}
\; \; \;
{\mathcal A}_{N \backslash G} =
{\mathcal A}_G \cap
A_{N \backslash G} .
$$
Note that $G/N, \; N \backslash G$ shall be denoted more precisely
by $G/(N, \pi), \; (N, \pi) \backslash G$ respectively, if there is
a possible confusion. Let $h_N$ be the Haar measure on $N$.
Let
$$
E_{G/N} = ( id \otimes h_N \pi ) \Delta
, \; \; \;
\; \; \;
E_{N \backslash G} = (h_N \pi \otimes id) \Delta .
$$
 Then   $E_{G/N}$
and $E_{N \backslash G}$ are projections of norm one
(completely positive and completely bounded conditional expectations)
from $A_G$ onto
$ A_{N \backslash G}$ and $ A_{G/N}$ respectively
(cf. \cite{Pod6} as well as Proposition 2.3 and Section 6 of \cite{W14}),
and
$$
{\mathcal A}_{G/N} = E_{G/N}( {\mathcal A}_G ) , \; \; \;
\text{and}
\; \; \;
{\mathcal A}_{N \backslash G}  =
E_{N \backslash G} ( {\mathcal A}_G ) .
$$
From this, we see that the *-subalgebras
$ {\mathcal A}_{N \backslash G}$ and $ {\mathcal A}_{G/N}$ are dense in
$A_{G/N}$ and $A_{N \backslash G}$ respectively.

Assume $N$ is a closed subgroup of an ordinary compact group $G$.
Let $\pi$ be the restriction morphism from $A_G : = C(G)$ to $A_N := C(N)$.
Let $C(G/N)$ and $C(N \backslash G)$ be continuous functions on
$G/N$ and $N \backslash G$ respectively.
Then one can verify that
$$
C(G/N) =
{A}_{G/N}= E_{G/N}(A_G),
$$
$$C(N \backslash G) =
{ A}_{N \backslash G}=
 E_{N \backslash G}(A_G).
$$
Therefore we will use the symbols
$C(G/N)$ and
${A}_{G/N}$
(resp. $C(N \backslash G)$ and ${ A}_{N \backslash G}
$; $C(G)$ and $A_G$)
interchangeably for all quantum groups.

\bgprop
\label{normal-subgroup}
Let $N$ be a quantum subgroup of a compact quantum group $G$.
Then the following conditions are equivalent:

{\rm (1)}  $ A_{N \backslash G}$ is a Woronowicz $C^*$-subalgebra of $A_G$.

{\rm (2)}  $ A_{G/N}$ is a Woronowicz $C^*$-subalgebra of $A_G$.

{\rm (3)}  $ A_{G/N} = A_{N \backslash G}$.

{\rm (4)}  For every irreducible representation $u^\lambda$ of $G$, either
$h_N \pi (u^\lambda) = I_{d_\lambda}$ or $h_N \pi (u^\lambda) = 0$,
where $h_N$ is the Haar measure on $N$, $d_\lambda$ is the dimension of
$u^\lambda$ and $I_{d_\lambda}$ is the $d_\lambda \times d_\lambda$
identity matrix.
\ndprop
\pf
We only need to show that
(1)$\Leftrightarrow$(4)$\Leftrightarrow$(3).
The proof of the implications
(2)$\Leftrightarrow$(4)$\Leftrightarrow$(3) is similar.

(3)$\Rightarrow$(4):
In general one has
$$
\Delta ({\mathcal A}_{N \backslash G})
\subseteq {\mathcal A}_{N \backslash G} \otimes {\mathcal A}_{G} , \; \; \;
\Delta( {\mathcal A}_{G/N})
\subseteq  {\mathcal A}_{G} \otimes {\mathcal A}_{G/N} .
$$
Letting ${\mathcal B} = {\mathcal A}_{N \backslash G} = {\mathcal A}_{G/N}$
then one has
$$
\Delta ({\mathcal B}) \subseteq
{\mathcal B} \otimes {\mathcal B} .
$$
For $\lambda \in \hat{G}$,
let $n_\lambda$ be the multiplicity of the trivial representation
of $N$ in the representation $\pi(u^\lambda)$.
We claim that either $n_\lambda = d_\lambda$ or
$n_\lambda = 0$.

Assume on the contrary that there is a $\lambda \in \hat{G}$ such that
$1 < n_\lambda < d_\lambda$.
Note that  in general
$$
E_{N \backslash G}(u^\lambda_{ij})
 =  (h_N \pi \otimes id) \Delta (u^\lambda_{ij})
 = \sum_k h_N \pi(u^\lambda_{ik}) u^\lambda_{kj},
$$
$$
E_{G/N}(u^\lambda_{ij})  = ( id \otimes h_N \pi ) \Delta(u^\lambda_{ij})
 = \sum_k  h_N \pi( u^\lambda_{kj}) u^\lambda_{ik}.
$$
Using unitary equivalence if necessary we choose
$u^\lambda_{ij}$ in such a way that
the $n_\lambda$ trivial representations of $N$ appear on the upper left
diagonal corner of $\pi(u^\lambda)$. Then
$$\displaystyle{
E_{N \backslash G}(u^\lambda_{ij}) =
                    \begin{cases}
                 u^\lambda_{ij} & \text{if \ $1 \leq i \leq n_\lambda , \;
                                             1 \leq j \leq d_\lambda  ,$}\\
                 0             & \text{if \ $n_\lambda < i \leq d_\lambda , \;
                                             1 \leq j \leq d_\lambda ,$}
                     \end{cases}}
$$

$$ \displaystyle{
E_{G/N}(u^\lambda_{ij}) =
                  \begin{cases}
                  u^\lambda_{ij} & \text{if \ $1 \leq i \leq d_\lambda , \;
                                             1 \leq j \leq n_\lambda  ,$} \\
                  0             & \text{if \ $1 \leq i \leq d_\lambda , \;
                                             n_\lambda < j \leq d_\lambda .$}
                  \end{cases}}
$$

Since
$ {\mathcal A}_{N \backslash G} = {\mathcal A}_{G/N} = {\mathcal B} $ and both
$E_{N \backslash G }$ and $E_{G/N}$
are projections from ${\mathcal A}_G$ onto ${\mathcal B}$, we have
$$
E_{N \backslash G} = E_{G/N}.
$$
Then for $n_\lambda < j \leq d_\lambda$,
$$
0 \neq u^\lambda_{ij} = E_{N \backslash G}(u^\lambda_{ij})
= E_{G/N}(u^\lambda_{ij})= 0,
$$
This is a contradiction.

(4)$\Rightarrow$(3):
Let $S(N)$ (or $S(N, \pi)$) be the subset of $\hat{G}$
consisting of those $\lambda$'s such that
$h_N \pi (u^\lambda)$ is $I_{d_\lambda}$. Then a straightforward calculation
using the fact that $E_{N \backslash G}$
and $E_{G/N}$ are 
projections of ${\mathcal A}_G$ onto
${\mathcal A}_{N \backslash G}$ and ${\mathcal A}_{G/N}$ respectively,
one gets
$$
 {\mathcal A}_{N \backslash G} = {\mathcal A}_{G/N}
= \bigoplus \{ {\mathbb C} u^\lambda_{ij} \;
| \; \lambda \in S(N), i, j = 1, \cdots, d_\lambda \} .
$$

(4)$\Rightarrow$(1):
Let $S(N)$ be defined as in the proof of (4)$\Rightarrow$(3).
It is clear that $A_{N \backslash G}$ is a Woronowicz $C^*$-subalgebra
of  $A_G$ and that
$$\{ u^\lambda_{ij} \; | \; \lambda \in S(N), i, j = 1, \cdots , d_\lambda
\}.
$$
is a Peter-Weyl basis of the dense $*$-subalgebra
${\mathcal A}_{N \backslash G}$ of ${A}_{N \backslash G}$.

(1)$\Rightarrow$(4):
Let $G_1 = N \backslash G$. Then by Woronowicz's Peter-Weyl theorem for
compact quantum groups,  every irreducible representation
$u^\lambda$ of $G$ is either an irreducible representation of $G_1$
or none of the coefficients $u^\lambda_{ij}$ is in ${\mathcal A}_{G_1}$.
That is
$$\displaystyle{
E_{N \backslash G} (u^\lambda_{ij})
                    = \begin{cases}
                    u^\lambda_{ij} & \text{if $\lambda \in \hat{G}_1$, } \\
                    0             & \text{if $
 \lambda \in \hat{G} \backslash \hat{G}_1 $}.
                    \end{cases}
               }
$$
By the definition of $E_{N \backslash G}$ and linear independence of
the $u^\lambda_{ij}$'s, this implies that
$$
h_N \pi (u^\lambda_{ik}) = \delta_{ik}, \; \; \lambda \in \hat{G}_1,
i, k = 1, \cdots , d_\lambda ,
$$
$$
h_N \pi (u^\lambda_{ik}) = 0, \; \;
 \lambda \in \hat{G} \backslash \hat{G}_1 .
$$
This completes the proof of the proposition.
\qed

\bgdf
\label{dfnormalsubgroup}
A quantum subgroup $N$ of a compact quantum group
$G$ is said to be {\bf normal} if it satisfies
the equivalent conditions of \propref{normal-subgroup}.
\nddf
\noindent
{\bf Remarks.}
(a) Condition (4) of \propref{normal-subgroup}
plays an important role in this paper.
It is a reformulation of the following condition
for a normal quantum subgroup $N$ that appears near the
end of Sect. 2 of \cite{W1}:
For every irreducible representation $v$ of $G$,
the multiplicity of the trivial representation of $N$ in the representation
$\pi(v)$  is either zero or the dimension of $v$.
From the proof of the proposition we see that the counit of
${\mathcal A}_{G/N}$ is equal to the restriction morphism
$\pi|_{{\mathcal A}_{G/N}}$.

(b) Note also that on p679 of \cite{W1} the following statement is found:
``In general, a right quotient
quantum group is different from the corresponding left
quotient quantum group.'' Though in the purely algebraic setting of
Hopf algebras, one needs to distinguish between left and right normal quantum
subgroups, as indicated in 1.5 of Parshall and Wang \cite{ParshallWang91a}
(see also \cite{Andrus95a,Schneid93a,Takeuchi94a}),
however, in view of \propref{normal-subgroup} above, this cannot happen
for normal quantum subgroups of compact quantum groups.
Moreover, using \lmref{full Woronowicz $C*$-algebras}-\lmref{reconstruct N from G/N} below, it can be
shown that the notion of normality defined in \cite{ParshallWang91a}
when applied to compact quantum groups is equivalent to our notion of normality.
As the main results of this paper do not depend on this equivalence, its proof and
other applications are in \cite{normal}.

(c) The notion of a normal quantum subgroup depends on the morphism $\pi$,
which gives the ``position'' of the quantum group $N$ in $G$.
If $(N, \pi_1)$ is another quantum subgroup of $G$ with surjection
$\pi_1: A_G \longrightarrow A_N$,  $(N, \pi_1)$ may not be normal
even if  $(N, \pi)$ is. This phenomenon already occurs in the
group situation. For example a
finite group can contain two isomorphic subgroups
with one normal but the other not.
\vv
{\bf Examples.} We show in (1) and (2) below that the identity group and the
full quantum group $G$ are both normal quantum subgroups of $G$
under natural embeddings.
These will be called the {\bf trivial normal quantum subgroups}.
See \secref{B_u-A_{aut}}-\secref{K_q-K_J} and
\cite{normal} for examples of non-trivial normal quantum subgroups.

(1) Let $N= \{ e \}$ be the one element identity group.
Let $\pi= \epsilon=$ counit of $A_G$
be the morphism from $A_G$ to $A_N$.
Then by the counital property, one has
$$A_{G/N} = \{ a \in A_G | (id \otimes \epsilon) \Delta (a) = a \otimes 1\}
= A_G.
\; \; \; $$
That is $(\{ e \}, \epsilon)$ is normal and
$G/(\{ e \}, \epsilon) =
 G$.

(2) Now let $N = G$ and let $\pi: A_G \rightarrow A_N$ be any isomorphism of
Woronowicz $C^*$-algebras \cite{W1}.
Let $h$ be the Haar measure on $G$ and $a \in A_{G/N}$.
Since $\pi$ is an isomorphism
and $(id \otimes \pi) \Delta (a) = a \otimes 1$,
one has $  \Delta(a) = (id \otimes \pi)^{-1} (a \otimes 1) = a \otimes 1$.
From the invariance of $h$ and  one has
$$
h(a) 1 = (1 \otimes h) \Delta(a) = a h(1) = a.
$$
Hence
$$A_{G/N} = \{ a \in A_G | (id \otimes \pi) \Delta (a) = a \otimes 1 \}
= {\mathbb C} 1.
\; \; \; $$
That is $(G, \pi)$ is normal and $G/(G, \pi) \cong
 \{ e \}$.

(3) We note that besides the embeddings in (2)
it is possible to construct examples of
compact quantum groups $G$ with non-normal proper embeddings of
$G$ into $G$. In fact this can happen for compact groups already.
\qed
\vv
\indent
 The following is a justification of the
above notion of normal quantum subgroups.

\bgprop
\label{justifydfnormalsubgroup}
Let $A=C(G)$ with $G$ a compact group. Let $N$ be a closed subgroup of $G$.
Let $\pi$ be the restriction map from $A$ to $A_N = C(N)$. Then $(N, \pi)$
is normal in the sense above if and only if $N$ is a normal subgroup of $G$
in the usual sense.
\ndprop
\pf
Under the Gelfand-Naimark correspondence which
associates to every commutative $C^*$-algebra
its spectrum, quotients of $G$
by (ordinary) closed normal subgroups $N$
correspond to Woronowicz $C^*$-subalgebras of $C(G)$, i.e.,
$$
G/N \; \; \; \text{corresponds to} \; \; \;
C(G/N),
$$
see 2.6 and 2.12 of \cite{W1}.
Since $A_{G/N}= C(G/N)$ for any closed subgroup $N$,
the proposition follows from \propref{normal-subgroup} above.
\qed
\vv
The following result gives a complete description of quantum
normal subgroups of the compact quantum group dual
of a discrete group $\Gamma$, whose proof
is straightforward using e.g.
\cite{Wor5} and \propref{normal-subgroup}.

\bgprop
\label{for-FD}
Let $A_G = C^*(\Gamma)$.
Let $N$ be a quantum subgroup of $G$ with surjection
$\pi: A_G \rightarrow A_N$.
Then $N$ is normal, $L:= \pi(\Gamma)$ is a discrete group
and $A_N = C^*(L)$.
Moreover, $A_{G/N} = C^*(K)$, where $K = \ker(\pi: \Gamma \rightarrow L)$.
\ndprop

To distinguish two different quantum subgroups, we include
the following result, which should be known to experts in
the theory of $C^*$-algebras.

\bgprop
\label{same}
Let $\pi_k: A \rightarrow A_k$ be surjections of unital
$C^*$-algebras with kernels $I_k$ ($k=1, 2$). Let $P_k$ be the
pure state space of $A_k$. Then the following
conditions are equivalent:

{\rm (1)}
$\{ \phi_1 \circ \pi_1 | \phi_1 \in P_1 \}
= \{ \phi_2 \circ \pi_2 | \phi_2 \in P_2 \}$ as subsets of pure states
of $A$.

{\rm (2)}
$I_1 = I_2$.

{\rm (3)}
There is an isomorphism $\alpha: A_1 \rightarrow A_2$ such that
$\pi_2 = \alpha \circ \pi_1$.
\ndprop
\pf
(1) $\Rightarrow$ (2): If $I_1 \neq I_2$, say, there is a
$x \in I_1 \backslash I_2$. Then there is a pure state $\phi$ of $A/{I_2}$
such that $\phi (\pi_2(x)) \neq 0$, where we identify
$A_2$ with $A/{I_2}$. But $\phi \pi_2$ is a pure state of $A/{I_1} \cong A_1$
according to assumption (1). Hence we must have $\phi \pi_2(x)=0$.
This is a contradiction.

(2) $\Rightarrow$ (3): Let $I=I_1=I_2$. Let $\pi$ be the
quotient map $A \rightarrow A/I$. Let $\tilde{\pi}_k$ be the
homomorphism from $A/I$ to $A_k$ such that
$\pi_k = \tilde{\pi}_k \pi$ ($k=1, 2$). Then $\tilde{\pi}_k$ are
isomorphisms. Put $\alpha = \tilde{\pi}_2 \circ {\tilde{\pi}_1}^{-1}$.
Then $\pi_2 = \alpha \circ \pi_1$.

(3) $\Rightarrow$ (1): This follows from
$P_1 = P_2 \circ \alpha$.
\qed
\vv
The following proposition is an easy consequence of
\propref{normal-subgroup}.

\bgprop
Let $(N_1, \pi_1)$
be a normal quantum subgroup of $G$. Let
$$\alpha: A_{N_1} \longrightarrow A_{N_2}$$
be an isomorphism of quantum groups. Then
$(N_2, \alpha \pi_1)$ is normal.
\ndprop
In view of the above discussions,
it is reasonable to have the following definition (cf. also remarks
after \propref{justifydfnormalsubgroup}).

\bgdf
\label{dfsame}
Two quantum subgroups $(\pi_1, H_1)$ and
$(\pi_2, H_2)$ of a quantum group $G$
are said to have {\bf the same imbedding} in $G$
if $\pi_1, \pi_2$ satisfy the equivalent conditions of
\propref{same}.  When this happens, we
denote  $( H_1 , \pi_1) = (H_2 , \pi_2)$.
\nddf

Geometrically speaking,
two quantum subgroups $( H_1 , \pi_1) $ and
$(H_2 , \pi_2)$ of a quantum group $G$
are said to have {\bf the same imbedding} in $G$
if their ``images'' in $G$ are the same.

\section{Simple Compact Quantum Groups}
\label{simple}

To avoid such difficulty such as the classification of finite groups
up to isomorphism in developing the theory of simple compact quantum groups,
we assume connectivity as a part of the postulates of the latter.
We use representation theory to define the notion of connectivity:

\bgdf
\label{connected}
We call a compact quantum group $G_A$ {\bf connected}
if for each non-trivial irreducible representation
$u^\lambda  \in \hat{G}_A$, the $C^*$-algebra $C^*(u^\lambda_{ij})$
generated by the coefficients of $u^\lambda$
is of infinite dimension.
\nddf

In virtue of (28.21) of \cite{HR2}, we have

\bgprop
\label{connected-prop}
Let $G_A$ be an ordinary
compact group (i.e. $A_G$ is commutative). Then $G_A$ is connected as a
topological space if and only if it is connected in the sense above.
\ndprop

\bgdf
\label{df-simple}
We call a compact
quantum group $G_A$ {\bf simple}
if it satisfies the following conditions {\rm (1)-(4):}

{\rm (1) }
The Woronowicz $C^*$-algebra $A_G$
is finitely generated;

{\rm (2) }
$G_A$ is connected;

{\rm (3) }
$G_A$ has no non-trivial
connected
normal quantum subgroups;

{\rm (4) }
$G_A$ has no non-trivial representations of dimension one.

A (simple) quantum group is called {\bf absolutely simple} if it
has no non-trivial normal quantum subgroups. Similarly a finite quantum group
is called {\bf simple} if it has no non-trivial normal quantum subgroups.
\nddf

Just as the notion of simple compact Lie groups excludes the torus groups,
the above notion of simple quantum groups excludes abelian compact quantum
groups in the sense of Woronowicz \cite{Wor5},
i.e.
quantum groups coming from group $C^*$-algebras $C^*(\Gamma)$ of discrete groups
$\Gamma$ (note that $C^*(\Gamma)$ is the algebra
of continuous functions on the torus ${\mathbb T}^n$
when $\Gamma$ is the discrete group ${\mathbb Z}^n$).
This is important because it is impossible to classify discrete
groups up to isomorphism.
However, we do not know if condition (4) in
\dfref{df-simple} (i.e., there is  no non-trivial group-like elements)
is superfluous, as is the case for simple compact Lie groups because of the
Weyl dimension formula.

As a justification of this definition, we have the following
proposition that shows that our notion of simple compact quantum groups
recovers {\em exactly} the ordinary notion of simple compact
Lie groups.

\bgprop
\label{justify-simple}
If $G_A$ is a simple compact quantum group with $A$ commutative, then
the set $G := \hat{A}$ of Gelfand characters is a
simple compact Lie group in the ordinary sense.
Conversely, every simple compact Lie group  in the ordinary sense
is of this form.
\ndprop

The proof \propref{justify-simple} follows immediately from Theorem 2.8 in \cite{W1}
and \propref{connected-prop} above.
We remark that although it is easy as above to prove the characterization of
the ordinary simple compact Lie group in terms of our notion of
simple compact quantum groups when $A_G$ is commutative,
it has been highly non-trivial to prove the analogous
characterization of ordinary differential manifolds
in terms of the axioms of non-commutative manifolds
which is finally achieved in the recent work of Connes \cite{Connes08mfd}
(see references therein for earlier, presumably un-successful,
attempts to at such a characterization).

Note that a simple compact Lie group is not a direct product of
proper connected subgroups. Also, a simple Lie group is not a
semi-direct product. Similarly, the following general results are true
for quantum groups (for proofs see \cite{normal}):

\bgprop
If $G_A$ is a simple compact
quantum group, then $A_G$ is not a tensor product, nor a crossed product
by a non-trivial discrete group.
\ndprop

To put in perspective the examples of simple compact quantum groups
to be studied later, we introduce some properties
for compact quantum groups.
First we recall that the {\bf representation ring} (also called the
{\bf fusion ring}) $R(G)$
of a compact quantum group $G$ is an ordered algebra over the
integers ${\mathbb Z}$ with positive cone
(or semiring, which is also a basis)
$R(G)_+: = \{ \chi_u \}$ consisting of characters
$\chi_u$ of irreducible representations $u \in \hat{G}$ of $G$,
and structure constants $c_{uv}^w \in {\mathbb N} \cup \{ 0 \}$
given by the rules
$$
\chi_{u} \chi_{v} = \sum_{w \in \hat{G}} c_{uv}^w \chi_{w},
$$
where the product $\chi_{u} \chi_{v}$ is taken in
the algebra ${\mathcal A}_G$.

\bgdf
\label{property-F-etc}
Let $G$ be a compact quantum group.
We say that $G$ has {\bf property $F$} if each
Woronowicz $C^*$-subalgebra of $A_G$ is of the form $A_{G/N}$ for some
normal quantum subgroup $N$ of $G$.
We say that $G$ has {\bf property $FD$}
if each quantum subgroup of $G$ is normal.

We say that $G$ is {\bf almost classical} if its representation ring $R(G)$
is order isomorphic to the representation ring of
a compact group.
\nddf

By \propref{justifydfnormalsubgroup}, a compact group trivially
has property $F$.
We will give in \secref{B_u-A_{aut}} and \secref{K_q-K_J} non-trivial
simple compact quantum groups that are almost classical and have property $F$.
Among compact quantum groups,
simple compact quantum groups that are almost classical or
have property $F$ are closest to ordinary simple compact Lie groups
in regard to noncommutative geometry.

By \propref{for-FD}, as the dual of discrete group $\Gamma$,
a compact quantum group of the form $C^*(\Gamma)$
has property $FD$. When $\Gamma$ is finite, $C^*(\Gamma)$
is equal to the dual of the function algebra $C(\Gamma)$.
This explains the term $FD$.

A compact quantum group $G$ is absolutely simple with property $F$
if and only if every non-trivial representation $v$ of $G$ is faithful, i.e.,
$C^*(v_{ij})= A_G$, cf. \cite{W14}.

By a theorem of Handelman \cite{Handel94a}, the representation
ring of a compact connected Lie group is a complete isomorphism invariant.
But this fails for compact quantum groups, since
the representation rings of a simple compact Lie group $K$
and its standard deformation $K_q$ are order isomorphic.

In \cite{Banica99c}, Banica uses the positive cone $R_+(G)$
of the representation ring $R(G)$ of a compact quantum group $G$ to define
what he calls an $R_+$ deformation.
This is closely related to almost classical quantum groups.

It is clear that a quantum quotient group $G/N$ of an almost classical quantum
group $G$ is also almost classical. But a quantum subgroup of an almost
classical quantum group need not be almost classical. For example,
the quantum permutation groups are almost classical
(cf. \cite{Banica7,W15} and remarks preceding \thmref{simpleA_{aut}}),
but according to of Bichon \cite{Bichon2}, their quantum subgroups
$A_2({\mathbb Z}/m{\mathbb Z})$
are not almost classical if $m \geq 3$ (see Corollary 2.7 and
the paragraph following Corollary 4.3 of \cite{Bichon2}).
However, for a compact quantum group with property $F$,
we have the following general result.

\bgth
Let $G$ be a compact quantum group with property $F$
Then its quantum subgroups and quotient groups $G/N$
(by normal quantum subgroups $N$) also have property $F$.
\ndth

As we will only use the definitions of quantum groups with
property $F$ (resp. property $FD$) but not the assertion in the theorem above,
the details for the proof of the theorem is included in
a separate paper \cite{normal}.
\vv
{\bf The main goals/problems} in the theory of simple compact quantum groups
are: (1) to construct and classify (up to isomorphism if possible)
simple compact quantum groups; (2) to construct and classify
irreducible representations of simple compact quantum groups;
(3) to analyze the structure of compact quantum
groups in terms of simple ones; and (4) to
develop applications of simple compact quantum groups
in other areas of mathematics and physics. For these purposes,
new techniques for compact quantum groups must be developed.

The above is a very difficulty program at the present.
Even problem (1) of the program above is daunting.
To obtain {\em clues} on the general problem (1), it is desirable
to find and solve easier parts of it. For this purpose, we propose
the following apparently easier problems.

\bgprob
(1) \label{prob-simple-property-F}
Construct and classify all simple compact quantum groups with property $F$ (up to isomorphism if possible).

(2) \label{prob-simple-almost-classical}
Construct and classify all simple compact quantum groups that are almost classical
(up to isomorphism if possible).
\ndprob

\bgprob
\label{prob-simple-property-FD}
Construct simple compact quantum groups with property $FD$.
\ndprob

Simple quantum groups in \probref{prob-simple-property-F} -- \probref{prob-simple-property-FD}
are most closest to groups known in mathematics. They should be easiest
classes to classify. Therefore they should play a fundamental role
in the main problems in the theory of simple compact quantum groups.

\section{Simplicity of
$B_u(Q)$ and $A_{aut}(B, \tau)$}
\label{B_u-A_{aut}}

To prove the main results in this  section and the next section,
we develop here two technical results,
which are of interest in their own right:
one on the correspondence between Hopf $*$-ideals and
Woronowicz $C^*$-ideals (\lmref{Hopf $*$-ideals vs. Woronowicz $C^*$-ideals});
the other on the reconstruction of a normal
quantum group from the identity in the quotient quantum group (\lmref{reconstruct N from G/N}).

We first recall the construction of compact quantum group
$B_u(Q)$ associated to a non-singular $n \times n$
complex scalar matrix $Q$ (cf. \cite{Banica1,W1,W5,W5'}).
The (noncommutative) $C^*$-algebra of functions on the quantum group
$B_u(Q)$ is generated by noncommutative coordinate functions
$u_{ij}$ ($i,j = 1, \cdots , n$) that are
subject to the following relations:
\bgeqq
    u^* u = I_n = u u^*, \; \; \;
    u^t Q u Q^{-1} = I_n = Q u Q^{-1} u^t,
\ndeqq
where $u=(u_{ij})_{i,j=1}^n$.
When $Q \bar{Q}$  is a scalar multiple $ c I_n$ of the identity
matrix $I_n$, the quantum group $B_u(Q)$ and the group $SU(2)$
have the same fusion rules for their irreducible representations,
as shown by Banica \cite{Banica1}, which
implies that $B_u(Q)$ is an almost classical
quantum group. Under the condition $Q \bar{Q} = \pm I_n$,
the isomorphism classification of $B_u(Q)$ is determined by the author
\cite{W17} using polar decomposition of $Q$ and eigenvalues
of $|Q|$ (see Theorem 2.4 in \cite{W17}).
For arbitrary $Q$, $B_u(Q)$ is a free product of its building blocks, involving
both $B_u(Q_l)$'s and $A_u(P_k)$'s
with $Q_l \bar{Q_l} $ being scalar matrices and
$P_k$ positive matrices (see Theorem 3.3 in \cite{W17}).
The precise definition of $A_u(Q)$ is recalled later in
the paragraphs before \propref{non-simpleA_u(Q)}.
For positive matrix $Q$,  $A_u(Q)$ is classified up to
isomorphism in terms of the eigenvalues of $Q$ (see Theorem 1.1 in \cite{W17});
and for a arbitrary non-singular matrix $Q$,
the general $A_u(Q)$ is a free products of $A_u(P_k)$'s
with positive matrices $P_k$ (see Theorems 3.1 in \cite{W17}).
In Bichon {\sl et al.} \cite{Bi-R-V06}, the same techniques in  \cite{W17}
were used to classify the unitary fiber functors of
the quantum groups $A_u(Q)$ and $B_u(Q)$ and their
ergodic actions with full multiplicity.
Note that for $n=1$, $B_u(Q)=C({\mathbb T})$ is the trivial
$1 \times 1$ unitary group. We will concentrate on the non-trivial
case $n \geq 2$. Note that the isomorphism class of $B_u(Q)$
depends on the normalized $Q$ only if $Q \bar{Q}$ is a scalar matrix
\cite{W17}.

\bgth
\label{simpleB_u(Q)}
Let $Q \in GL(n, {\mathbb C})$ be such
that $Q \bar{Q} = \pm I_n$. Then $B_u(Q)$ is
an almost classical
simple compact quantum group with property $F$.
In fact it has only one normal subgroup of order $2$.
\ndth
\pf
As noted above, the quantum group $B_u(Q)$ is almost classical
because its representation ring is order isomorphic to the representation ring
of the compact Lie group $SU(2)$ \cite{Banica1}.
More precisely, according to \cite{Banica1} irreducible representations
of the quantum group $B_u(Q)$ can be
parametrized by $r_k$ ($k = 0, 1, 2, \cdots$)
with $r_0$ trivial and $r_1 = (u_{ij})_{ij=1}^n$, so that the fusion rules
for their tensor product representations
(i.e., decomposition into irreducible representations)
read
$$
r_k \otimes r_l = r_{|k-l|} \oplus r_{|k-l| + 2 } \oplus \cdots
\oplus r_{k+l -2} \oplus r_{k+l}, \; \; k, l \geq 0.
$$

We show that the quantum group $B_u(Q)$ is connected.
If $k=2m$ is even ($m >0$), then let $r_l=r_k$ in the above
tensor product decomposition and do the same for the
irreducible constituents repeatedly, one sees that the algebra
$C^*(r_{2m})$ generated
by the coefficients of the representation $r_{2m}$
contains the coefficients of $r_{2s}$ for all $s$.
Hence $r_{2m}$ generates an infinite dimensional algebra:
$$C^*(r_{2m})= C^*(\{ r_{2s} | s \geq 0 \} ).$$
If  $k=2m+1$ is odd ($m \geq 0$), then let $r_l=r_k$ in the above
tensor product decomposition, one sees that the representation
$r_2$ appears
therein.
Apply the decomposition to
$r_{2m+1} \otimes r_2$, one sees that $r_1 = (u_{ij})$
appears therein. Hence the algebra generated by the coefficients
of $r_{2m+1}$ is the same as the algebra generated by those
of $r_1=(u_{ij})$.
We conclude from this analysis that {\em there is only one non-trivial
 Woronowicz $C^*$-subalgebra
in $B_u(Q)$, the one $C^*(r_{2m})$
generated by coefficients of $r_{2m}$,}
which is obviously infinite dimensional as noted above,
where $m$ is any nonzero positive number.
In particular, the quantum group $B_u(Q)$ is connected.

For rest of the proof,
we show that the quantum group $B_u(Q)$ has only one
normal quantum subgroup, although it has many quantum subgroups.

Note that the coordinate functions $v_{ij}$
of the matrix group $N = \{ I_n,  - I_n \}$
satisfy the defining relations of $B_u(Q)$, hence
there is a surjection $\pi$ from the $C^*$-algebra $B_u(Q)$
to the $C^*$-algebra $A_N$
of functions on $N$ such that
$$
\pi(u_{ij}) = v_{ij}, \; \; \; i,j = 1, 2, \cdots, n.
$$
It is clear that $\pi$ is a morphism of quantum groups, hence
$(N, \pi)$ is a quantum subgroup of the quantum group
$B_u(Q)$.

We show that $(N, \pi)$ is actually a normal quantum subgroup.
To see this, it suffices by \propref{normal-subgroup}
to show that
$$
\pi(r_{2m}) = d_{2m} \cdot v_0, \; \; \;
\pi(r_{2m+1}) = d_{2m+1} \cdot v_1,
$$
where $d_{2m}$ and $d_{2m+1}$ are dimensions of the representations
$r_{2m}$ and $r_{2m+1}$ respectively,
$v_0$ and $v_1$ are the trivial representation and  the
non-trivial irreducible representation of $N$
respectively ($v_1(\pm I_n) = \pm 1$).
By the definition of $\pi$ and $v_1$
the assertion is clearly true for $m=0$.
In general, suppose the assertion is true for $m$.
Then $\pi(r_{2m +1}) \otimes \pi(r_1)$ is a multiple of
$v_0$ since  $v_1^2 = v_0$. From the decomposition of
 $r_{2m +1} \otimes r_1$,
we get
$$
\pi(r_{2m +1}) \otimes \pi(r_1) =
\pi(r_{2m}) \oplus \pi(r_{2m +2}).
$$
Hence $\pi(r_{2(m+1)}) = \pi(r_{2m +2})$ is a scalar
multiple of $v_0$.
Similarly,  from
$$
\pi(r_{2m + 2}) \otimes \pi(r_1) =
\pi(r_{2m + 1}) \oplus \pi(r_{2m + 3}),
$$
we see that
$ \pi(r_{2(m+1) + 1})= \pi(r_{2m + 3}) $
is a multiple of $v_1$. Since $v_0$ and
$v_1$ are one dimensional representations,
the multiples we obtained above must be
$d_{2m+2}$ and $d_{2m+3}$ respectively.
That is $(N, \pi)$ is normal and
$$
A_{G/N} = C^*(r_{2})= C^*(\{ r_{2s} | s \geq 0 \} ),
$$
where for simplicity of notation,
the symbol $G$ in $G/N$ refers to the quantum group $G_{B_u(Q)}$.
The above also shows that this quantum group has property $F$.

We have to show that $B_u(Q)$ has no other normal quantum subgroups,
which will imply that it has no connected normal quantum
subgroups and is therefore a simple quantum group.

Let $(N_1, \pi_1)$ be a non-trivial normal quantum subgroup of $B_u(Q)$.
We show that $(N_1, \pi_1) = (N, \pi)$ in the sense of \dfref{dfsame}, which will finish the proof
of the theorem.
Since 
$N_1 \neq 1$, by \dfref{dfsame} and \propref{normal-subgroup}
there exists an irreducible representation $v$
of the quantum group $B_u(Q)$ such that $\pi_1(v)$ is not a scalar
and therefore $E_{G/{N_1}}(v)=0$.
Hence by the proof of \propref{normal-subgroup} and Woronowicz's
Peter-Weyl theorem \cite{Wor5},
$A_{G/{N_1}} = E_{G/{N_1}}(A_G) \neq A_G$.

Similarly, we claim that
$A_{G/{N_1}} \neq {\mathbb C} 1$, where $1$ is the unit of $A_G$.
To prove this, we need three lemmas.
It is instructive to compare the second lemma
(\lmref{Hopf $*$-ideals vs. Woronowicz $C^*$-ideals})
with the ideal theory for $C^*$-algebras.

\bglm
\label{full Woronowicz $C*$-algebras}
Let $B_1$ and $B_2$ be Woronowicz $C^*$-algebras with canonical dense
Hopf $*$-algebras of ``representative functions'' ${\mathcal B_1}$ and ${\mathcal B_2}$ respectively.
Assume $B_2$ is full and $\psi : B_1 \rightarrow B_2$ is a morphism of Woronowicz $C^*$-algebras such that
the induced morphism $\hat{\psi} : {\mathcal B_1} \rightarrow {\mathcal B_2}$ is an isomorphism.
Then $B_1$ is full and $\psi$ is also an isomorphism.
\ndlm
\noindent
{\bf Remark.} The above is false if the roles of $B_1$ and $B_2$ are exchanged, as seen by taking
$B_1= C^*(F_2)$ and $B_2= C^*_r(F_2)$.
\vv
{\em Proof of \lmref{full Woronowicz $C*$-algebras}}.
Since ${\mathcal B_1}$ is dense in $B_1$, it suffices to show that $||\hat{\psi}(a)|| = ||a||$ for
$a \in {\mathcal B_1}$.

Since $\psi$ is a morphism of $C^*$-algebras,
we have $||{\psi}(a)|| \leq ||a||$ and therefore the first
inequality
$$
||\hat{\psi}(a)|| = ||{\psi}(a)|| \leq ||a|| \,.
$$

Since $B_2$ is full, the norm on ${\mathcal B_2}$ is the universal
$C^*$-norm (see \cite{W3}):
$$
||\hat{\psi}(a)|| = \sup \{ ||\pi(\hat{\psi}(a))||: \pi \;
\text{is a *-representation of $\mathcal B_2$} \} \,.
$$
Taking $\pi= \hat{\psi}^{-1}$ in the above, we obtain the second inequality
$$
||\hat{\psi}(a)|| \geq ||\hat{\psi}^{-1}(\hat{\psi}(a))|| =||a|| \, .
$$
Combining the first the second inequalities finishes
the proof of \lmref{full Woronowicz $C*$-algebras}.

\bglm
\label{Hopf $*$-ideals vs. Woronowicz $C^*$-ideals}
{\rm (Hopf $*$-ideals vs. Woronowicz $C^*$-ideals)}

{\rm (1) }
Let $G$ be a compact quantum group. Let ${\mathcal I}$ be a
Hopf $*$-ideal of ${\mathcal A}_G$.
Then the norm closure $\overline{\mathcal I}$ in the 
$C^*$-algebra $A_G$ is a Woronowicz $C^*$-ideal and
$A_G/\overline{\mathcal I}$ is a full Woronowicz $C^*$-algebra.
The Hopf $*$-algebra ${\mathcal A}_G/{\mathcal I}$ admits a universal
$C^*$-norm and its completion under this norm is a
Woronowicz $C^*$-algebra isomorphic to
$A_G/\overline{\mathcal I}$.

{\rm (2) }
The map $f({\mathcal I})=\overline{\mathcal I}$ is a bijection
from the set of Hopf $*$-ideals $\{ {\mathcal I} \}$
of ${\mathcal A}_G$ onto the set of Woronowicz
$C^*$-ideals $\{ I \}$ of $A_G$ such that $A_G / I$ is full.
The inverse $g$ of $f$ is given by
$g(I) = I \cap {\mathcal A}_G$.
\ndlm
\noindent
{\bf Remarks.}
(a) Note that (2) and the last part of (1) in
the lemma above are false if the Woronowicz $C^*$-algebra
$A_G$ or $A_G/I$ is not full, as is shown by the following example.
Let $A_G=C^*(F_2)$ be the group $C^*$-algebra of the
free group $F_2$ on two generators.
Let $I$ be the kernel of the canonical map
$\pi: C^*(F_2) \rightarrow C^*_r(F_2)$
where $C^*_r(F_2)$ is  the reduced group $C^*$-algebra of $F_2$.
Then $I \cap {\mathcal A}_G=0$ but $ I \neq \overline{0}$.

(b) This lemma strengthens the philosophy in \cite{W3} that the ``pathology'' associated with
the ideals between $0$ and the kernel of the morphism from the full Woronowicz $C^*$-algebra to
reduced one such as $\pi: C^*(F_2) \rightarrow C^*_r(F_2)$ is not (quantum) group theoretical,
but purely functional analytical, and $C^*(F_2)$ and $C^*_r(F_2)$ should be viewed as the same
quantum group because the same dense Hopf $*$-subalgebra
that completely determines the quantum group can be recovered from either the full or the reduced
algebra.
Similarly, for a general compact quantum group $G$,
the totality of {\em (quantum) group theoretic} information
is encoded in the {\em purely algebraic} object ${\mathcal A}_G$,
any other (Hopf) algebra should be viewed as defining the same quantum
group as ${\mathcal A}_G$ so long as ${\mathcal A}_G$ can be recovered from it.
The advantage of working with the category of full $C^*$-algebras or the purely algebraic
objects ${\mathcal A}_G$ is that morphisms can be easily defined for them,
whereas it is not even possible to define a morphism from the one element group
to the quantum group associated with the reduced algebra $C^*_r(F_2)$
if is viewed as a different quantum group than the one
associated with the full algebra $C^*(F_2)$.
\vv
{\em Proof of \lmref{Hopf $*$-ideals vs. Woronowicz $C^*$-ideals}}.

Let ${\mathcal I}$ be as in (1).
Let $\pi_1: A_G \rightarrow A_G/\overline{\mathcal I}$ be the quotient map.
Since ${\mathcal I}$ is a Hopf *-ideal, we have in particular (see Sweedler \cite{Sweedler})
$$
\Delta({\mathcal I})
\subset
{\mathcal A}_G \otimes {\mathcal I} + {\mathcal I} \otimes {\mathcal A}_G
\subset 
\ker(\pi_1 \otimes \pi_1) .
$$
Therefore
$
\Delta(\overline{\mathcal I}) \subset  \ker(\pi_1 \otimes \pi_1) .
$
That is, $\overline{\mathcal I}$ is a Woronowicz $C^*$-ideal
and $A_G/\overline{\mathcal I}$ is a Woronowicz
$C^*$-algebra (see 2.9-2.11 in \cite{W1}).
Denote $B_1 = A_G/\overline{\mathcal I}$ and let $\hat{\pi}_1$
be the induced morphism of the canonical dense Hopf-*-subalgebras
$\hat{\pi}_1: {\mathcal A}_G \rightarrow {\mathcal B_1}$.

We claim that $\ker{\hat{\pi}_1} = {\mathcal I}$ and
$\hat{\psi}_0: {\mathcal A}_G/{\mathcal I} \longrightarrow {\mathcal B}_1$,
$\hat{\psi}_0: [a] \mapsto \pi_1 (a)$
is an isomorphism, where
$[a] \in {\mathcal A}_G/{\mathcal I}$, $a \in {\mathcal A}_G$.

By \cite{Wor5,W1}, ${\mathcal A}_G$ is generated as an algebra
by the coefficients $u^\lambda_{ij}$ of irreducible unitary
corepresentations $u^\lambda$ of Hopf $*$-algebra ${\mathcal A}_G$.
The images $[u^\lambda_{ij}]$ of $u^\lambda_{ij}$ in the quotient Hopf $*$-algebra
${\mathcal A}_G/{\mathcal I}$ give rise to
unitary corepresentation of ${\mathcal A}_G/{\mathcal I}$,
and generate it as an {\em algebra} (not just as a *-algebra).
Therefore ${\mathcal A}_G/{\mathcal I}$
is a compact quantum group algebra (CQG algebra) in the sense
of Dijkhuizen and Koornwinder \cite{DijKwd1} (See also
\cite{KlimykSchmudgen} and \cite{Wor6,W3})--a more appropriate name
for compact quantum group (CQG) algebra might be Woronowicz $*$-algebra
(or compact Hopf $*$-algebra),
since the quantum group $C^*$-algebra of a compact quantum group $G$
is the $C^*$-algebra $C^*(G)$ dual to $C(G)$ according
to \cite{PW}.

Let ${\mathcal B }_2 = {\mathcal A}_G/{\mathcal I}$
and let $B_2$ be the closure of ${\mathcal B }_2$
in the universal $C^*$-norm. Then $B_2$ is a Woronowicz
$C^*$-algebra. As the norm on $A_G$ is universal, the composition
$$
{\mathcal A}_G
\longrightarrow {\mathcal A}_G/{\mathcal I}
\longrightarrow B_2
$$
is bounded and extends to a morphism of
Woronowicz $C^*$-algebras $\rho: A_G \rightarrow B_2$.
Since ${\mathcal I} \subset \ker(\rho)$, we have
$\overline{\mathcal I} \subset \ker(\rho)$ and
$\rho$ factors through $B_1 = A_G/\overline{\mathcal I}$ via a
$C^*$-algebra morphism $\psi$:
$$
A_G \stackrel{{\pi}_1} \longrightarrow B_1
\stackrel{\psi} \longrightarrow B_2 \, , \;  \rho = {\psi} {\pi}_1 .
$$
It is clear that $\rho(a) = [a]$ for $a \in {\mathcal A}_G$
and from this it can
be checked that $\hat{\psi}$ and $\hat{\psi}_0$ are inverse morphisms, where
$\hat{\psi}: {\mathcal B}_1 \longrightarrow  {\mathcal B}_2$ is the restriction of
$\psi$ to the dense Hopf $*$-subalgebra ${\mathcal B}_1$ and $B_2$.
Hence $\hat{\psi}_0$ is an isomorphism as claimed.

From $\hat{\rho} = \hat{\psi} \hat{\pi}_1$ (since $\rho = {\psi} {\pi}_1$),
 it is easy to see that ${\psi}$ is a morphism of Woronowicz
$C^*$-algebras (see 2.3 in \cite{W1}). Since $\hat{\psi} = \hat{\psi}_0^{-1}$ is
an isomorphism and $B_2$ is full, by \lmref{full Woronowicz $C*$-algebras},
$ B_1 $ is full and $\psi$ is itself an isomorphism from $ B_1 $ to $B_2$.
(We note in passing that since $A_G/\ker(\rho) \cong B_2$,
we have $\overline{\mathcal I} = \ker(\rho)$.)
This proves part (1) of the lemma.

To prove part (2) of the lemma, let ${\mathcal I}$ be as in (2) and
$B_1= {A}_G/\overline{\mathcal I}$.
Then by (1) above and \cite{W1}, $B_1$ is a Woronowicz $C^*$-algebra.
Let ${\mathcal B}_1$ be the canonical dense Hopf $*$-algebra of $B_1$ and let
$\hat{\pi}_1: {\mathcal A}_G \rightarrow {\mathcal B}_1$ be
the morphism associated with the quotient morphism $\pi_1$.
 Then clearly
$${\mathcal I} \subset \overline{\mathcal I} \cap {\mathcal A}_G
=gf({\mathcal I}) .$$
 Conversely if $x \in \overline{\mathcal I} \cap {\mathcal A}_G  $,
then $x \in \ker(\hat{\pi}_1)
= {\mathcal I}$.
Hence $gf({\mathcal I}) = {\mathcal I}$.

Next let ${I}$ be as in (2).
We show that $fg(I) = I$. Let $B_2 = A_G/I$ -- this is not the same $B_2$ as in (1) above.
Let $\pi_2$ be the quotient morphism from  $A_G$ onto $B_2$
(compare with $\rho$ above).
Define ${\mathcal I} = g(I)=I \cap {\mathcal A}_G$.
We need to show that $\overline{\mathcal I} = I$. The idea of proof
is the same as that of the last part in (1).

Using the morphism
$\hat{\pi}_2: {\mathcal A}_G \rightarrow {\mathcal B}_2$
of dense Hopf *-algebras associated with $\pi_2$, we see that ${\mathcal I}= \ker(\hat{\pi}_2)$.
Hence ${\mathcal I}$ is a Hopf $*$-ideal in ${\mathcal A}_G$ and
${\mathcal A}_G/{\mathcal I}$ is isomorphic to ${\mathcal B}_2$ under the
natural map induced from $\hat{\pi}_2$, and by (1) above, $B_1 := A_G/\overline{\mathcal I}$ is a
Woronowicz $C^*$-algebra.
Since $\overline{\mathcal I} \subset I$, the morphism $\pi_2$ factors
through $B_1$ via a morphism $\psi$ of Woronowicz
$C^*$-algebras:
$$
A_G \stackrel{{\pi}_1} \longrightarrow B_1
\stackrel{{\psi}} \longrightarrow B_2 \, , \;  \pi_2 = {\psi} {\pi}_1 .
$$
Besides being isomorphic to ${\mathcal B}_2$,
${\mathcal A}_G/{\mathcal I}$ is also isomorphic to ${\mathcal B}_1$ (under the
morphism $\hat{\psi}_0$) according to the proof of (1) earlier.
Hence the restriction $\hat{\psi}$ of $\psi$ to the dense Hopf $*$-algebras
is an isomorphism from ${\mathcal B}_1$ to ${\mathcal B}_2$.
Since $B_2$ is full, by \lmref{full Woronowicz $C*$-algebras},
 $\psi$ itself is an isomorphism,
which means that $I = \overline{\mathcal I}$ (and $B_1=B_2$).
This completes the proof of
\lmref{Hopf $*$-ideals vs. Woronowicz $C^*$-ideals}.

\bglm
\label{reconstruct N from G/N}
{\rm (Reconstruct $N$ from $G/N$)}

Let $(N, \pi)$ be a normal quantum subgroup of a compact quantum
group $G$. Let $\hat{\pi}$ be the associated morphism from
${\mathcal A}_G$ to ${\mathcal A}_N$. Then,
$$
\ker(\hat{\pi}) =
{\mathcal A}^+_{G/{N}} {\mathcal A}_G =
{\mathcal A}_G {\mathcal A}^+_{G/{N}}
= {\mathcal A}_G {\mathcal A}^+_{G/{N}} {\mathcal A}_G,
$$
where  ${\mathcal H}^+$ denotes the
augmentation ideal (i.e. kernel of the counit)
for any Hopf algebra ${\mathcal H}$.
\ndlm
\noindent
{\bf Remarks.}
(a)
In the notation of Schneider  \cite{Schneid93a},
the result above can be restated as follows:
The map $\Phi$ is the left inverse of $\Psi$,
where  $\Psi(\ker(\hat{\pi})) := {\mathcal A}_{G/N}$ and
$\Phi({\mathcal A}_{G/N}):= {\mathcal A}_G {\mathcal A}_{G/N}^+ $.
In the language of Andruskiewitsch {\sl et al}
\cite{Andrus95a}, the result above implies that the sequence
$$
1 \longrightarrow
N  \longrightarrow
G \longrightarrow
G / N \longrightarrow
1,
$$
or the sequence
$$
0 \longrightarrow
{\mathcal A}_{G/{N}} \longrightarrow
{\mathcal A}_G \longrightarrow
{\mathcal A}_N \longrightarrow
0,
$$
is exact. It is instructive to compare this with the purely algebraic
situation in Parshall and Wang \cite{ParshallWang91a},
where for a given a normal quantum subgroup in the sense there,
the existence of an exact sequence is not known
and the uniqueness does not hold in general (cf. 1.6 and 6.3 loc. cit.).
Note that the notion of exact sequence of quantum groups in Schneider
\cite{Schneid93a} is equivalent to that in Andruskiewitsch {\sl et al} \cite{Andrus95a}
under certain faithful (co)flat conditions.
Though a Hopf algebra is not faithfully flat over its Hopf
subalgebras if it is not commutative or cocommutative
(see Schauenburg \cite{Schauenburg2000}), we
have the following
\vv
{\bf Conjecture:}
Let $G$ be a compact quantum group. Then the Hopf algebra
$A_G$ (resp. ${\mathcal A}_G$) is faithfully flat over its Hopf subalgebras.

Similarly, $A_G$ (resp. ${\mathcal A}_G$) is faithfully coflat over
$A_G/I$ (resp. ${\mathcal A}_N/{\mathcal I}$) for every Woronowicz $C^*$-ideal $I$
(resp. Hopf *-ideal ${\mathcal I}$).
\\

(b)
It can be shown using \lemref{reconstruct N from G/N}
and Schneider \cite{Schneid93a} that the notion of normal
quantum groups in this paper (or in \cite{Wang,W1}) and the one
in Parshall and Wang \cite{ParshallWang91a} are equivalent for compact quantum groups.
For more details, see \cite{normal}
\vv
{\em Proof of \lmref{reconstruct N from G/N}}.
The proof is an adaption of the ones in 16.0.2 of
Sweedler \cite{Sweedler} and (4.21) of Childs \cite{Childs} for finite
dimensional Hopf algebras to infinite dimensional ones considered here.
We sketch the main steps here for convenience of the reader.

It suffices to prove $
\ker(\hat{\pi}) =
{\mathcal A}^+_{G/{N}} {\mathcal A}_G$. The other equality
$
\ker(\hat{\pi}) =
{\mathcal A}_G {\mathcal A}^+_{G/{N}}
$
is proved similarly. From these it follows that
$
\ker(\hat{\pi}) = {\mathcal A}_G {\mathcal A}^+_{G/{N}} {\mathcal A}_G.
$

Consider the right ${\mathcal A}_N$-comodule structures on
${\mathcal A}_N$ and ${\mathcal A}_G$ given respectively
by
$$
\Delta_N:
{\mathcal A}_N  \rightarrow {\mathcal A}_N \otimes {\mathcal A}_N ,
 \; \; \text{and} \; \; \;
(id \otimes \hat{\pi})\Delta_G:
{\mathcal A}_G  \rightarrow {\mathcal A}_G \otimes {\mathcal A}_N ,
$$
where $\Delta_N$ and $\Delta_G$ are respectively the coproducts
of the Hopf algebras ${\mathcal A}_N$ and ${\mathcal A}_G$.
Since ${\mathcal A}_N$ is cosemisimple by the fundamental work of
Woronowicz \cite{Wor5} (see remarks in 2.2 of \cite{W1}), it follows from Theorem 3.1.5
of \cite{Dascalescu} that every
${\mathcal A}_N$-comodule is projective.
Furthermore, one checks that the surjection
$\hat{\pi}: {\mathcal A}_G \rightarrow {\mathcal A}_N$
is a morphism of
${\mathcal A}_N$-comodules. Hence $\hat{\pi}$ has a comodule
splitting
$s:  {\mathcal A}_N \rightarrow {\mathcal A}_G $
with $\hat{\pi} s = id_{{\mathcal A}_N}$.

Let $x \in {\mathcal A}^+_{G/{N}}$.
By remark (a) following \dfref{dfnormalsubgroup}, $\hat{\pi}(x) = 0$.
Hence  $ {\mathcal A}^+_{G/{N}} \subset \ker(\hat{\pi})$ and therefore
  $ {\mathcal A}^+_{G/{N}} {\mathcal A}_G \subset \ker(\hat{\pi})$.

Define a linear map $\phi$ on ${\mathcal A}_G$ by
$\phi = (s \hat{\pi}) * S = m (s \hat{\pi} \otimes S) \Delta_G$,
where $m$ and $S$ are respectively the multiplication map
and antipode of ${\mathcal A}_G$.
Then using the coassociativity of $\Delta_G$ and
$\hat{\pi} s = id_{{\mathcal A}_N}$ along with the antipodal property
of $S$, one verifies that $\phi({\mathcal A}_G) \subset {\mathcal A}_{G/{N}}$.
Since $\ker(\hat{\pi}) \subset Im(id- s \hat{\pi})$, to
show  $  \ker(\hat{\pi})  \subset {\mathcal A}^+_{G/{N}} {\mathcal A}_G$,
it suffices to show that
$Im(id- s \hat{\pi}) \subset {\mathcal A}^+_{G/{N}} {\mathcal A}_G$.
Since $(\epsilon-id) \phi({\mathcal A}_G) \subset {\mathcal A}^+_{G/N}$,
the later follows from the identity
$$
id- s \hat{\pi} =
(\epsilon-id) \phi * id
=m((\epsilon-id) \phi \otimes id)\Delta_G,
$$
which one verifies using basic properties of the convolution product
along with $\epsilon \phi = \epsilon$ and the splitting property of $s$.
This proves  \lmref{reconstruct N from G/N}.
\vv
\indent
Now we finish the proof of \thmref{simpleB_u(Q)}.
 If $A_{G/{N_1}} = {\mathbb C} 1$,
we would have ${\mathcal A}_{G/{N_1}} = {\mathbb C} 1$ and
${\mathcal A}^+_{G/{N_1}} =0$.
Let $\hat{\pi}_1$ be the morphism of Hopf algebras from
${\mathcal A}_G$ to ${\mathcal A}_{N_1}$ associated with $\pi_1$.
Then by \lmref{reconstruct N from G/N},
$$
\ker(\hat{\pi}_1) =
{\mathcal A}^+_{G/{N_1}} {\mathcal A}_G = 0.
$$
Since $\ker(\hat{\pi}_1)$ is dense in $\ker(\pi_1)$ by
\lmref{Hopf $*$-ideals vs. Woronowicz $C^*$-ideals}, we would have
$\ker(\pi_1)=0$. This contradicts the assumption that $N_1$ is a
non-trivial quantum subgroup of $G$ and therefore $A_{G/{N_1}}
\neq {\mathbb C} 1$.

Then $A_{G/{N_1}}$ has to be the only
non-trivial Woronowicz $C^*$-subalgebra of $B_u(Q)$, i.e.
$A_{G/{N_1}}= C^*(\{ r_{2m} | m \geq 0 \} )$
as noted near the beginning of the proof of \thmref{simpleB_u(Q)}.
We infer from \propref{normal-subgroup} that
$\pi_1(r_{2m})$ is a multiple of the trivial representation
of $N_1$ for any $m$. From
$$
\pi_1(r_{1}) \otimes \pi_1(r_1) =
\pi_1(r_{0}) \oplus \pi_1(r_{2}),
$$
we see that
$
\pi_1(r_{1}) \otimes \pi_1(r_1)
$
is a multiple of the trivial representation of $N_1$.
That is
$$
\sum_{ijkl} e_{ij} \otimes e_{kl} \otimes {\tilde u}_{ij} {\tilde u}_{kl}
= I_n \otimes I_n \otimes 1,
$$
where ${\tilde u}_{ij}$ are the
$(i,j)$-entries of $ \pi_1(r_1)$ and $e_{ij}$ are matrix units.
Hence
$$
{\tilde u}_{ij} {\tilde u}_{kl} = 0, \; \;
\text{when} \; \; i \neq j, \; \;
\text{or}
\; \;  k \neq l;
$$
$$
{\tilde u}_{ii} {\tilde u}_{ll} = 1, \; \;
\text{for all} \; \;
i, l.
$$
Therefore ${\tilde u}_{ij}=0$ for $i \neq j$ and
$A_{N_1}$ 
is commutative. That is, $N_1$ is an ordinary compact group.
Now it is clear that
${\tilde u}_{ii} = {\tilde u}_{ll}={\tilde u}_{ll}^{-1}$ for all $i,l$,
which we denote by $a$.
Since $A_{N_1}$ 
is generated by $a$ and $N_1$ is
non-trivial, we conclude that $N_1$ is a group
of order $2$.
The map $\alpha$ from
$A_N$ 
to $A_{N_1}$  
defined
by $\alpha(v_{ij})= \tilde u_{ij}$
is clearly an isomorphism such that $\pi_1 = \alpha \pi$.
Hence 
$(N_1, \pi_1) = (N, \pi)$ by \dfref{dfsame}.

For an example of a
non-normal quantum subgroup $(H, \theta)$ of $B_u(Q)$,
take a two-elements group $H = \{ I_n, V \}$, where
$$V =
\displaystyle{
\left[
\begin{array}{cc}
             -1 & 0    \\
             0  & I_{n-1}
\end{array}
\right],
}
$$
and $\theta(u_{ij}) = w_{ij}$, the coordinate functions on $H$.
\qed
\vv

Let us also recall the construction of the quantum groups
$A_u(Q)$ closely related to $B_u(Q)$ \cite{W1,W5,W5'}.
For every non-singular matrix $Q$,
the quantum group $A_u(Q)$ is defined in terms of generators
$u_{ij}$ ($i,j =1,  \cdots  n$), and relations:
\bgeqq
    u^* u = I_n = u u^*, \; \; \;
    u^t Q {\bar u} Q^{-1} = I_n = Q {\bar u } Q^{-1} u^t.
\ndeqq

According to Banica \cite{Banica2}, when $Q>0$,
the irreducible representations of
the quantum group $A_u(Q)$ are parameterized by
the free monoid ${\mathbb N} \ast {\mathbb N}$ with generators
$\alpha$ and $\beta$ and anti-multiplicative involution $\bar{\alpha}= \beta$
(the neutral element is $e$ with $\bar{e} = e$).
The classes of $u$ and $\bar{u}$ are $r_\alpha$ and $r_\beta$ respectively.
Moreover, for each pair of irreducible representations
$r_x$ and $r_y$ ($x , y \in {\mathbb N} \ast {\mathbb N}$),
one has the following direct sum decomposition (fusion) rules:
$$r_x \otimes r_y = \sum_{x = ag, \bar{g} b = y}  r_{ab}. $$
In \cite{W17},  the special $A_u(Q)$'s with $Q>0$ are classified
up to isomorphism and the general $A_u(Q)$'s with arbitrary
$Q$ are shown to be free products of the special $A_u(Q)$'s.
The following result was observed by Bichon through private communication
(the proof given below was developed by the author):

\bgprop
\label{non-simpleA_u(Q)}
{\em The quantum groups $A_u(Q)$ are not simple for
any $Q \in GL(n, {\mathbb C})$}.
\ndprop

\pf
To prove this, we first introduce the following notion.
A quantum subgroup $(N, \pi)$ of a compact quantum group $G$ is said to be {\bf in the center}
of $G$ if
$$
(\pi \otimes id ) \Delta = (\pi \otimes id )\sigma  \Delta \;,
$$
where $\sigma(a_1 \otimes a_2) = a_2 \otimes a_1$, $a_1, a_2 \in A_G$, and $\Delta$ is the coproduct of $A_G$.

Assume $(N, \pi)$ is in the center of $G$. Then using the definitions of
$ A_{G/N}$ and $A_{N \backslash G}$ in \secref{notion-of-normal},
it is straightforward to verify that $ A_{G/N} = A_{N \backslash G}$.
By \propref{normal-subgroup}(3), $(N, \pi)$ normal in $G$.
Namely, a quantum subgroup that is in the center of $G$ is always normal, just as in the classical case.

Let ${\mathbb T}$ be the one dimensional (connected) torus group and
$t \in C({\mathbb T})$ the function such that $t^*t =1= t t^*$. Then
$C({\mathbb T})$ is generated by $t$ as a $C^*$-algebra: $C({\mathbb T})= C^*(t)$.
Define the morphism $\pi: A_u(Q) \rightarrow C({\mathbb T})$ by
 $\pi(u_{ij}) = \delta_{ij} t$ (note the special case $A_u(Q)=C({\mathbb T})$ when $n=1$).
Then it is routine to verify that the connected group $(T, \pi)$ is in the center
of the quantum group $A_u(Q)$ (not viewed as an algebra) in the sense above
and is therefore a normal subgroup therein. Hence
$A_u(Q)$ is not simple.
\qed
\vv
\indent
We remark that although $A_u(Q)$ is not simple, for  $n \geq 2$ and $Q>0$
it is very close to being normal, satisfying
most of the axioms of a simple compact quantum group:
its function algebra is finitely generated; it is connected; and its non-trivial
irreducible representations are all of dimensional greater than one (see Wang \cite{W17}
for a computation of the dimension of its irreducible representations based on
Banica \cite{Banica2}). In particular
following problems should not be hard:

\bgprob
\label{simpleA_u(Q)}
(1) Study further the structure of $A_u(Q)$ for
positive matrices $Q \in GL(n, {\mathbb C})$ and $n \geq 2$.
Determine all of their simple quotient quantum groups. Alternatively,

(2)
\label{prob-almost-classical}
Construct simple compact quantum groups that are not almost classical.
\ndprob

A solution of part (1) of the above problem should also
give a solution to part (2) and  provide
the first examples of simple compact quantum groups that
are not almost classical because of the highly non-commutativity of the
representation ring of $A_u(Q)$
(note that all the simple quantum groups known so far are almost classical).
It is worth noting that the determination of all simple quotient quantum
groups of $A_u(Q)$ in the above problem is easier than the
determination of all of their simple quantum subgroups,
the latter being tantamount to finding all simple quantum groups
because every compact matrix quantum group is a quantum subgroup of
an appropriate $A_u(Q)$.
These remarks also indicate that that $A_u(Q)$ should play an important role
in the theory of simple compact quantum groups.

Next we consider the quantum automorphism group $A_{aut}(B, tr)$ of
a finite dimensional $C^*$-algebra $B$ endowed with
a tracial functional $tr$ (cf. \cite{Banica7,W15}).
This quantum group is defined to be the universal object
in the category of compact quantum transformation groups
of $B$ that leave $tr$ invariant.
Note that the presence of a tracial functional $tr$
is necessary for the existence of the universal object when $B$ is
non-commutative (see Theorem 6.1 of \cite{W15}).
For an arbitrary finite dimensional $C^*$-algebra $B$,
the $C^*$-algebra $A_{aut}(B, tr)$
is described explicitly in \cite{W15} in terms of generators and relations.
When $B=C(X_n)$ is the commutative $C^*$-algebra of functions on the space
$X_n$ of $n$ points, the quantum automorphism group $A_{aut}(B)= A_{aut}(X_n)$
(also called the {\bf quantum permutation group} on $n$ letters) exists without
the presence of a (tracial) functional and
its description in terms of generators
and relations is surprisingly simple. The $C^*$-algebra
$A_{aut}(X_n)$ is generated by self-adjoint
projections $a_{ij}$ such that each row and
column of the matrix $(a_{ij})_{i,j=1}^n$ adds up to 1. That is,
\bgeqq
\label{r1}  a_{ij}^2 = a_{ij} = a_{ij}^*,
\; \; \; i,j = 1, \cdots, n, \\
\label{r2}  \sum_{j = 1}^{n} a_{ij} = 1 ,
\; \; \; i = 1, \cdots, n, \\
\label{r3}  \sum_{i = 1}^{n} a_{ij} = 1 ,
\; \; \; j = 1, \cdots, n.
\ndeqq
For more general finite dimensional $C^*$-algebras $B$, the description
of $A_{aut}(B, tr)$ in terms of generators and relations is more complicated.
We refer the reader to \cite{W15} for details.

Assume $tr$ is the canonical trace $\tau$ on $B$
(see p772 of \cite{Banica7} or  \secref{Introd} for the definition).
Then $A_{aut}(B, \tau)$ is an ordinary permutation group
when the dimension of $B$ is less than or equal to $3$.
However, when the dimension of $B$ is greater than or equal to
$4$, $A_{aut}(B, \tau)$ is a non-trivial
(noncommutative and noncocommutative) compact quantum group
with an infinite dimensional function algebra \cite{W15,W14},
and as Banica  \cite{Banica7} showed, the algebra of symmetries
of the fundamental representation of this quantum group is
isomorphic to the infinite dimensional Temply-Lieb algebras $TL(n)$ and
the representation ring of $A_{aut}(B, \tau)$
is isomorphic to that of $SO(3)$. Hence
$A_{aut}(B, \tau)$ is almost classical for all $B$.
It is easy to see that for $B=C(X_n)$,
the canonical trace $\tau$ is equal to the
unique $S_n$-invariant state on $B$, where
$S_n$ acts on $X_n$ by permutation.
Hence by remark (2) following Theorem 3.1 of \cite{W15},
$A_{aut}(B, \tau)$ is the same as the quantum permutation
group $A_{aut}(X_n)$.

We refer the reader to \cite{Banica7,W15,W14}
for more on these quantum groups and
\cite{Bichon1,Bichon2,Bichon3} for interesting related results.
Note that the description in \cite{Banica7} is not exactly
as that in \cite{W15} but equivalent to it.
We now prove

\bgth
\label{simpleA_{aut}}
Let $B$ be a finite dimensional $C^*$-algebra with $\dim(B) \geq 4$.
Then $A_{aut}(B, \tau)$ is an almost classical, absolutely simple compact quantum group
with property $F$.
\ndth

\pf
The argument is similar to the one in \thmref{simpleB_u(Q)}. By Banica \cite{Banica7},
the complete set of mutually inequivalent irreducible representations
of the quantum group $A_{aut}(B, \tau)$
can be parametrized by $r_k$ ($k \geq 0$, $r_0$ being the trivial
one dimensional representation). Under this parametrization the fusion rules
of its irreducible representations are the same as those of $SO(3)$ and therefore
it is almost classical:
$$
r_k \otimes r_l = r_{|k-l|} \oplus r_{|k-l| + 1 } \oplus \cdots
\oplus r_{k+l -1} \oplus r_{k+l}, \; \; k, l \geq 0.
$$

We claim that there are only two Woronowicz
$C^*$-subalgebras in
$A_{aut}(B, \tau)$, namely ${\mathbb C}1$ and $A_{aut}(B, \tau)$.

Let $A_1 \neq {\mathbb C}1$ be a
Woronowicz $C^*$-subalgebra of $A_{aut}(B, \tau)$.
Let $v$ be a non-trivial irreducible representation of the compact
quantum group of $A_1$. Then $v=r_k$ for some
$k \neq 0$ and
$$
r_k \otimes r_k = r_{0} \oplus r_{ 1 } \oplus r_{ 2 } \oplus \cdots
\oplus r_{2k -1} \oplus r_{2k}, \; \; k, l \geq 0.
$$
Hence the coefficients of each of the representations
$r_1$, $r_2$, ..., $r_{2k}$ are in $A_1$.
Similarly, from the decomposition of $r_{2k} \otimes r_{2k}$,
we see that
the coefficients of each of the representations
$r_1$, $r_2$, ..., $r_{4k}$ are in $A_1$.
Inductively,
the coefficients of each of the representations
$r_1$, $r_2$, ..., $r_{2^m k}$ are in $A_1$ ($m>0$).
Hence $A_1 = A_{aut}(B, \tau)$.

Let $(\pi, N)$ be a normal quantum subgroup
of $G=A_{aut}(B, \tau)$ different from the trivial one-element
subgroup.
Then there is a non-trivial irreducible representation
$u^\lambda =(u^\lambda_{ij})$ such that $\pi(u^\lambda) $
is not a multiple of the trivial representation.
Using the same argument as in the proof of
\thmref{simpleB_u(Q)} we have
$A_{G/{N}} = E_{G/{N}}(A_G) \neq A_G$.
Therefore we must have $A_{G/N} = {\mathbb C} 1$.
Then the argument near the end of the proof of \thmref{simpleB_u(Q)}
(i.e. the paragraph that follows
the proof of \lmref{reconstruct N from G/N})
shows that $\ker(\pi)=0$. That is,  $N$ is the same quantum group as $G$.
\qed
\vv
\indent
\thmref{simpleA_{aut}} applies in particular to
quantum permutation groups $A_{aut}(X_n)$ when $n \geq 4$.
As Manin (private communication in July, 2002)
pointed out to the author (private communication), the reason that
these quantum groups are connected could be that there are so many
more quantum symmetries that the originally $n!$ isolated
permutations are connected together by them. Note however that
their function algebras are generated by orthogonal projections
$a_{ij}$, so these quantum groups are also disconnected, as
observed by Bichon \cite{Bichon2}.
It would be interesting to find a satisfactory explanation of
this paradox.

The proofs of the main results of this section do not
need explicit description (models) of representations of the quantum groups
$B_u(Q)$ and $A_{aut}(B, \tau)$ and $A_u(Q)$. Only the structures of
their representation rings (i.e. fusion rules) are used.
However, explicit constructions of models
of irreducible representations of Lie groups are fundamental
and have important applications in other branches of mathematics and physics.
Moreover, just as the construction and classification of the
representations of simple compact Lie groups is intimately intertwined with
the classification of simple compact Lie groups,
the same might hold true for simple compact quantum groups.
In view of these, we believe an appropriate answer to the following
problem should be important in the theory of compact quantum groups
in general and the theory of simple compact quantum groups in particular.
(Note that the model for the fundamental representation of the
quantum group $A_u(Q)$ is used in \cite{W14} to construct
ergodic actions on various von Neumann factors.)

\bgprob
Construct explicit models of the irreducible
representations of the following quantum groups:
$A_u(Q)$ for $Q>0$;
$B_u(Q)$ for $Q \bar{Q} = \pm I_n$;
the quantum automorphism group $A_{aut}(B, \tau)$ of
a finite dimensional $C^*$-algebra $B$ endowed with
the canonical trace $\tau$. Relate the results to the
theory simple compact quantum groups if possible.
\ndprob

\section{Simplicity of $K_q$, $K_q^u$ and $K_J$
}
\label{K_q-K_J}

The compact real forms $K_q$ of
Drinfeld-Jimbo quantum groups
and their twists $K_q^u$ are studied in \cite{Soib1} and
\cite{Lev1} respectively.
See also \cite{LS1} for a summary of \cite{Soib1,Lev1}
and \cite{KorogodskiSoib1} for more detailed treatment.
Motivated by these works, Rieffel constructs in \cite{Rf8} a
deformation $K_J$ of compact Lie group $K$ which contains a torus $T$ and
raises the question whether $K^u_q$
can be obtained as a {\em strict deformation quantization} of $K_q$.
This question is answered in the affirmative by the author in \cite{W4}.
The purpose of this section is to show that the quantum groups $K_q$, $K^u_q$
and $K_J$ are simple in the sense of this paper, provided that the compact
Lie group $K$ is simple.

We first recall the notation of \cite{Soib1,Lev1,LS1}.
Let $G$ be a connected and simply connected
simple complex Lie group with Lie algebra $ \mathfrak{g} $.
Fix a triangular decomposition
$\mathfrak{g} = \mathfrak{n}_{-} \oplus \mathfrak{h} \oplus \mathfrak{n}_{+}$,
 together with the corresponding decomposition
 $\Delta = \Delta _{+} \cup \Delta _{-}$
of the root system and a fixed basis $\{ \alpha _i \}_{i=1}^{n}$
for $\Delta _+$.
For each linear functional $\lambda $ on $\mathfrak{h}$,
$H_\lambda$ denotes the element in $\mathfrak{h}$ corresponding to
$\lambda$ under the isomorphism $\mathfrak{h} \cong \mathfrak{h}^*$ determined
by the Killing form $( \; , \; )$ on $\mathfrak{g}$.
Note that if the reader keeps the context in mind,
the symbols $\alpha$ and $\lambda$ used in this context should not
cause confusion with the same symbols used in this paper for other purposes.
Let $\{ X_\alpha \}_{\alpha \in \Delta } \cup \{ H_i \}_{i=1}^n$ be
a Weyl basis of $\mathfrak{g}$, where $H_i = H_{\alpha _i}$.
This determines a Cartan involution $\omega _0$ on $\mathfrak{g}$ with
$\omega _0 ( X_\alpha ) = -X_{- \alpha }$,
$\omega _0 ( H_i ) = - H_i$. Let $\mathfrak{k}$ be the compact real form of
$\mathfrak{g}$ defined as the fixed points
of $\omega _0$ and $K$ the associated
compact real form of $G$.
Put ${\mathfrak h}_{\mathbb R} = \oplus_{i=1}^n {\mathbb R} H_i$,
${\mathfrak t} = i {\mathfrak h}_{\mathbb R}$
and $T = \exp({\mathfrak t })$, the later being
the associated maximal torus of $K$.

Let
$q=e^{h/4}$
($h \in {\mathbb R} \backslash \{ 0 \}$).
For $n, k\in {\mathbb N}$, $n\ge k$, define
\begin{eqnarray*}
[n]_{q} &=&\frac{q^n -q^{-n}}{q -q^{-1}},\\
{}\left[{n\atop k}\right]_{q}&=&\frac{[n]_{q}[n-1]_{q} \ldots {}[n-k+1]_{q}}
{[k]_{q}[k-1]_{q}\ldots [1]_{q}}.
\end{eqnarray*}

The quantized universal enveloping algebra $U_q( \mathfrak{g} )$
\cite{Dr1,Jimb1} is the complex associative algebra
with generators $X_i^{{}\pm{}}$, $K_i^{{}\pm 1}$\/
($i=1, \cdots, n$) and defining relations:
\bgeqq
& & K_i K_i^{- 1} = 1 =  K_i^{ -1} K_i,\; \;  K_i K_j = K_j K_i  , \\
& & K_i X_j^{{}\pm{}} K_i^{-1} =
q^{{}\pm (\alpha_i, \alpha_j)} X_j^{{}\pm{}}  , \\
& & {}[X_i^+ , X_j^-] = \delta_{ij} \frac{K^2_i - K_i^{-2}}{q-q^{-1}}  , \\
& & \sum_{k=0}^{1-a_{ij}}
 (-1)^k \left[ {1-a_{ij}\atop k}\right ]_{q_i} \;
(X_i^{{}\pm{}})^k X_j^{{}\pm{}}(X_i^{{}\pm{}})^{1-a_{ij}-k}= 0\;,\;i\ne j,
\ndeqq
where $q_i = q^{(\alpha_i , \alpha_i)}$.

On $U_q(\mathfrak{g})$\/ there is a Hopf algebra structure with
coproduct
$$ 
\Delta(K_i^{{}\pm 1}) = K_i^{{}\pm 1}\otimes K_i^{{}\pm 1} ,\; \; \;
\Delta(X_i^{\pm}) =  X_i^{\pm} \otimes K_i + K_i^{-1} \otimes X_i^{\pm} ,
$$ 
 and counit and antipode respectively
$$ 
\varepsilon (X_i^{\pm}) = 0 , \; \; \;
\varepsilon (K_i^{\pm 1}) = 1 , \;  \; \;
S (X_i^{\pm}) = -q_i^{\pm 1} X_i^{\pm}, \; \; \;
S (K_i^{\pm 1}) = K_i^{\mp 1}.
$$ 
Under the *-structure defined by
$$(X_i^{\pm})^* = X_i^{\mp}, \; \; \; K_i^* = K_i,$$
$U_q(\mathfrak{g})$ is a Hopf *-algebra.

Let $u = \sum_{k,l} c_{kl} H_k \otimes H_l
\in \wedge^2 \mathfrak{h}_{\mathbb{R}}$.
Then it can be shown (cf. \cite{KorogodskiSoib1}) that
the following defines a new coproduct on $U_q (\mathfrak{g})$,
$$\Delta_u (\xi) = \exp(-ihu/2) \Delta (\xi) \exp(ihu/2),$$
 where $X \in U_q(\mathfrak{g})$ and $\Delta$ is the original coproduct
 on $U_q(\mathfrak{g})$. The new Hopf *-algebra so obtained is denoted by
 $U_{q,u}(\mathfrak{g})$.

The function algebra ${\mathcal A}_{K_q}$ of the
compact quantum group $K_q$ is defined
to be the subalgebra of the dual algebra $\Uq^*$
consisting of matrix elements of finite dimensional representations $\rho$
of $\Uq$ such that eigenvalues of the endomorphisms $\rho(K_i)$ are positive.
The function algebra ${\mathcal A}_{K^u_q}$ of the compact quantum group
$K^u_q$ is defined to be the subalgebra of the dual algebra
 $U_{q,u}(\mathfrak{g})^*$ that has the same elements as
${\mathcal A}_{K_q}$, as well as the same $*$-structure,
while the product of its elements is defined using $\Delta_u$ instead
of $\Delta$.

For each (algebraically) dominant integral weight $\lambda \in  P_+$ of
$(\mathfrak{ g, h})$, define matrix elements $C^\lambda _{\mu, i; \nu ,j }$
of the highest weight  $\Uq$ module $(L(\lambda), \rho_\lambda)$ as follows.
Let $\{ v_\nu^{(i)} \}$ be an orthonormal weight basis for the
unitary $\Uq$ module $L(\lambda)$.
Then $C^\lambda _{\mu , i; \nu ,j }$ is defined by
$$
C^\lambda _{\mu , i; \nu ,j }(X) =
  < \rho_\lambda (X) v_\nu^{(j)},  v_\mu^{(i)}>,$$
where $X \in \Uq$ and $<\; , \;>$ is the inner product on $L(\lambda)$.
The $C^\lambda _{\mu, i; \nu ,j }$'s is a linear (Peter-Weyl)
basis of both ${\mathcal A}_{K_q}$ and
${\mathcal A}_{K^u_q}$ when $\lambda$ ranges through
the set $P_+$ of dominant integral weights of $(\mathfrak{ g, h})$.

\bgth
\label{simpleK_q}
Let $K$ be a connected and simply connected simple compact Lie group.
Then for each each $q$,  $K_q$ is an almost classical simple compact quantum group with property $F$.
\ndth
\pf
First we recall that representations of $K$ and $K_q$ are in one
to one correspondence via deformation
and the decompositions of tensor products of
irreducible representations are not altered under deformation
(see Lusztig \cite{Lus88a}
and Rosso \cite{Rosso88a} or Chari-Pressley \cite{ChariPressley}).
 From this it follows immediately that $K_q$  is almost classical.

Let $\xi$ be the 
map that associates each irreducible representation $v$
of $K$ an irreducible representation $\xi(v)$ of $K_q$ in this
correspondence. This map defines an isomorphism of {\em vector spaces} from
${\mathcal A}_{K}$ to  ${\mathcal A}_{K_q}$, which we also denote by $\xi$.
It follows from this that $K_q$ is connected and has no
non-trivial representations of dimension one.
Comparing decompositions of tensor products of representations
of $K$ and $K_q$ we see that the $\xi$ maps bijectively the set of
Hopf subalgebras of ${\mathcal A}_{K}$ onto the set of
Hopf subalgebras of ${\mathcal A}_{K_q}$.

Let $\rho_q$
be the quotient morphism from
$A_{K_q}$ to the abelianization $A^{ab}_{K_q}$, which is by definition
the quotient of $A_{K_q}$ by the closed two sided ideal of
$A_{K_q}$ generated by commutators $[a, b]$, $ a, b \in A_{K_q}$.
According to \cite{W1},  $A^{ab}_{K_q}$ is the algebra of continuous
functions on the maximal compact subgroup $\hat{A}_{K_q}$ of
$K_q$ and $\rho_q$
gives rise to the embedding of the quantum groups from $\hat{A}_{K_q}$
to $K_q$. It is shown in \cite{Soib1} that the maximal compact subgroup
$\hat{A}_{K_q}$ is isomorphic to the maximal torus $T$ of $K$.

The associated morphism $\hat{\rho}_q$ from
${\mathcal A}_{K_q}$ to ${\mathcal A}_{T}$ is given by
$$
\hat{\rho}_q(C^\lambda_{\mu, i; \nu, j})(t) =
\delta_{ij} \delta_{\mu \nu} e^{2 \pi \mu(x)},
$$
where $t = exp(x) \in T$, $x \in \mathfrak{t}= i \mathfrak{h}_{\mathbb R}$
(see p438 of \cite{ChariPressley}, but $\sqrt{-1}$ should not appear in
the formula there). It is clear that one has the same formula as above
for restriction morphism ${\rho}$
from ${ A}_{K}$ to ${ A}_{T}$:
$$
\hat{\rho}( \xi^{-1}(C^\lambda_{\mu, i; \nu, j}))(t) =
\delta_{ij} \delta_{\mu \nu} e^{2 \pi \mu(x)},
\; \; i.e., \; \;
\hat{\rho} = \hat{\rho}_q \circ \xi.
$$

Let $N \subset K$ be a normal subgroup of $K$ with surjections
$\pi: A_{K} \rightarrow A_{N}$ and
$\hat{\pi}: {\mathcal A}_{K} \rightarrow {\mathcal A}_{N}$.
Then $N$ is a finite subgroup of $T$ and
$A_{N}= {\mathcal A}_{N}$ is a finite dimensional Hopf algebra.
It is clear that $\pi =  \rho_N \circ \rho$,
where $\rho_N$ is the restriction morphism from $A_T$ to $A_N$.
Define
$$
\pi_q: A_{K_q} \longrightarrow A_N,
\; \; \text{by} \; \;
\pi_q := \rho_N \circ \rho_q.
$$
We claim that $(N, \pi_q)$ is a normal subgroup of $K_q$.
This follows immediately from the following identities,
which one can easily verify using
$\hat{\rho} = \hat{\rho}_q \circ \xi$
and
$\hat{\pi} = \hat{\pi}_q \circ \xi$.
$$
{\mathcal A}_{K_q/N} = \xi({\mathcal A}_{K/N}),
\; \; i.e.,
\; \;
$$
$$
\{ a \in {\mathcal A}_{K_q} | (id \otimes \pi_q) \Delta (a) = a \otimes 1\}
=
\xi (\{ a \in {\mathcal A}_K | (id \otimes \pi) \Delta (a) = a \otimes 1\});
$$
$$
{\mathcal A}_{N \backslash K_q} = \xi({\mathcal A}_{N \backslash K}),
\; \; i.e.,
\; \;
$$
$$
\{ a \in {\mathcal A}_{K_q} | (\pi_q \otimes id) \Delta (a) = 1 \otimes a \}
=
\xi (\{ a \in {\mathcal A}_K | (\pi \otimes id) \Delta (a) = 1 \otimes a \})
.$$
That is, every normal subgroup $N$ of $K$ gives rise to a normal subgroup
$(N, \pi_q)$ of $K_q$ in the manner above.

Conversely, let $(N', \pi')$ be a quantum normal subgroup of $K_q$.
Then
$
{\mathcal A}_{K_q/N'}
$
is a Hopf subalgebra of ${\mathcal A}_{K_q}$.
Since every Hopf subalgebra of ${\mathcal A}_K$ is of the
form ${\mathcal A}_{K/N}$ for some normal subgroup
$N$ of $K$ (cf. \cite{W1}),
by the correspondence between Hopf subalgebras of  ${\mathcal A}_K$
and those of  ${\mathcal A}_{K_q}$  noted near the beginning of the proof
we have
$$
{\mathcal A}_{K_q/N'} =
\xi({\mathcal A}_{K/N} ) =
{\mathcal A}_{K_q/N}
$$
for some normal subgroup $N$ of $K$. By
\lmref{reconstruct N from G/N},
we have
$\ker(\hat{\pi})= \ker(\hat{\pi}')$.
That is $(N', \pi')$ and $(N, \pi_q)$ are the
same quantum subgroup of $K_q$
(cf. \dfref{dfsame} and
\lmref{Hopf $*$-ideals vs. Woronowicz $C^*$-ideals}).
Since normal subgroups $N$ of $K$ are finite, we conclude from
the above that $K_q$ has no non-trivial connected quantum normal
subgroups.
\qed
\vv
Examining the proof of \thmref{simpleK_q}, we formulate the following
general result on the invariance of simplicity of compact quantum groups
under deformation, which will be used to prove
the simplicity of $K_q^u$ and $K_J$.

Let $G$ be an almost classical simple compact quantum group with property $F$
and $(H, \rho)$ a quantum subgroup.
Assume all normal quantum subgroups of $G$ are quantum subgroups of $H$.
Let $G_v$ be a family of compact quantum groups (``deformation'' of $G$)
indexed by a subset $\{ v \}$ of a vector space that includes the origin.
Suppose the family $G_v$ satisfies the following conditions:

C1. $G_0 = G$.

C2. There is an isomorphism $\xi$ of vector spaces from
${\mathcal A}_{G}$ to ${\mathcal A}_{G_v} $.

C3. The coproduct is unchanged under deformation, i.e.,
$$
\Delta_v(\xi(a)) = (\xi \otimes \xi) \Delta(a)
\; \;
\text{for}
\; \;
a \in {\mathcal A}_{G}.
$$

C4. For any pair irreducible representations $u^{\lambda_1} $
and $u^{\lambda_2}$ of $G$, if
$$
u^{\lambda_1} \otimes u^{\lambda_2} \cong
u^{\gamma_1} \oplus u^{\gamma_2} \oplus \cdots \oplus u^{\gamma_l}
$$
is a decomposition of $u^{\lambda_1} \otimes u^{\lambda_2} $
into direct sum of irreducible subrepresentations $u^{\gamma_j}$
($j=1, 2, \cdots, l$),  then
$$
\xi(u^{\lambda_1}) \otimes \xi(u^{\lambda_2}) \cong
\xi(u^{\gamma_1})  \oplus \xi(u^{\gamma_2}) \oplus
\cdots \oplus \xi(u^{\gamma_l})
$$
is a decomposition of $\xi(u^{\lambda_1}) \otimes \xi(u^{\lambda_2})$
into direct sum of irreducible representations, where for instance
$\xi(u^{\lambda_1})$ denotes the representation of $G_v$
whose coefficients are images of coefficients of $u^{\lambda_1}$.

C5.  The quantum subgroup $H$ is undeformed. The latter means that
there is a morphism $\rho_v$ of quantum groups from $H$ to $G_v$ such that
$$
\rho_v(\xi(a)) = \rho(a) \; \; \text{for } \; \; a \in {\mathcal A}_G.
$$

Under the assumptions above, we have the following result.
The proof is the same as that of \thmref{simpleK_q}
($H$ corresponds to $T$ in \thmref{simpleK_q}).

\bgth
\label{simple-invariance}
For each $v \in \{ v \} $,  $G_v$ is an almost classical simple quantum group with property $F$.
\ndth
\noindent
{\bf Remarks.}
(a)
Condition C4 above is not the same as the requirement that
$$
\xi(u^{\lambda_1} \otimes u^{\lambda_2}) =
\xi(u^{\lambda_1}) \otimes \xi(u^{\lambda_2}).
$$
The latter requirement together with conditions (2) and (3) imply that
$\xi$ is an isomorphism of quantum from $G_v$ to $G$, which is not the
case for the quantum groups under consideration here.

(b)
We believe similar results on invariance of simplicity under deformation
hold true without the property $F$ assumption on $G$.
But at the moment we do not know of any simple compact quantum groups
that do not satisfy this property, though there many non-simple quantum
groups without this property.
\vv
Next we recall the construction in \cite{Rf8,W4}.
Let $A=A_G$ be a compact quantum group with coproduct
$\Delta$. Suppose that the quantum group $G$ has a toral subgroup
$(T, \rho)$--to obtain non-trivial deformation we assume that $T$
has rank no less than $2$.
For any element $t$ in $T$, denote by $E_t$ the corresponding
evaluation functional on $C(T)$. Assume that $\eta$
is a continuous  homomorphism from a vector space Lie group
${\mathbb{R} }^n$ to $T$, where $n$ is allowed to
be different from the dimension of $T$.
Define an action $\alpha$ of
${\mathbb{R} }^d : =
{\mathbb{R} }^n \times {\mathbb{R} }^n$ on the $C^*$-algebra
$A$ as follows:
$$\alpha_{(s,v)} = l_{\eta(s)} r_{\eta(v)}, $$
where
$$l_{\eta(s)} = (E_{\eta{(-s)}} \rho \otimes id) \Delta, \; \; \;
  r_{\eta(v)} = (id \otimes E_{\eta(v)} \rho ) \Delta.$$
For any skew-symmetric operator $S$ on
${\mathbb{R} }^n$, one may apply Rieffel's quantization procedure
\cite{Rf6} for the action $\alpha$ above to
obtain a deformed $C^*$-algebra $A_J$ whose product is denoted $\times_J$,
where $J = S \oplus (-S)$.
The family $A_{h J}$ ($h \in {\mathbb R}$) is a strict deformation
quantization of $A$ (see Chapter 9 of \cite{Rf6}).
In \cite{W4} the following result is obtained.

\bgth
\label{rieffel-deform}
The deformation $A_J$ is a compact quantum group containing
$T$ as a (quantum) subgroup; $A_J$ is a compact matrix
quantum group if and only if $A$ is.
\ndth

We denote by ${G_J}$ the quantum group for $A_J$. When $G$ is a compact Lie
group, the construction $G_J$ above is the same as Rieffel's
construction \cite{Rf8}. By 5.2 of \cite{Rf8}, $G_J$ is an almost classical
compact quantum group if $G$ is a compact Lie group.

Combining \thmref{rieffel-deform} and \thmref{simple-invariance},
we obtain

\bgth
\label{simpleK_J}
Let $K$ be a simple compact Lie group with a toral subgroup $T$.
Then $K_J$ of Rieffel \cite{Rf8} is an almost classical simple
compact quantum group with property $F$.
\ndth

We note that unlike in \thmref{simpleK_q}, in the result above we do not
need to assume $K$ to be simply connected. This is because
${\mathcal A}_{K_q}$ is defined using irreducible representations of
$U_q( \mathfrak{g} )$ associated with all dominant integral weights $P_+$ of
$({\mathfrak g}, {\mathfrak h})$, so that
${\mathcal A}_{K_q}$ becomes the algebra of representative functions on
a simply connected $K$ when $q \rightarrow 1$.
One could also start with a non-simply connected
$K$ in \thmref{simpleK_q} too, but then one needs to modify the definition of
the quantum algebra ${\mathcal A}_{K_q}$ by
using irreducible representations of $U_q( \mathfrak{g} )$ associated with
analytically dominant integral weights only.
This newly defined ${\mathcal A}_{K_q}$ is a Hopf subalgebra of the
Hopf algebra defined originally.
It is clear from the proof of \thmref{simpleK_q}
that its conclusion remains valid for this newly defined $K_q$.

Finally we consider $K_q^u$.
To avoid confusion with the Killing form, we now use
$s \oplus v$, instead of $(s,v)$ used above,  to
denote an element of ${\mathbb{R} }^d =
 {\mathbb{R} }^n \times {\mathbb{R} }^n$.
In the present setting, the space ${\mathbb{R}}^n$ is
$ {\mathfrak{h}}_{\mathbb R}$,
with inner product $<\; , \; > = (\; , \;)$,
where $(\; , \;)$ is the Killing form of $ \mathfrak{g} $  restricted
 to $ {\mathfrak{h}}_{\mathbb R} $.
We will also use  $<\; , \;>$  to denote the inner product on
$ {\mathfrak{h}}_{\mathbb R} \oplus  {\mathfrak{h}}_{\mathbb R} $.
Noting that the compact abelian group $T$ is also a subgroup of
both $K_q$ and $K_q^u$ (see \cite{Soib1,Lev1}).
The map $\eta$ there in this case is
 defined by $\eta (s) = \exp( 2 \pi i s)$. We can define as above
an action of ${\mathbb R}^d$ on $A_{K_q}$ by
$$\alpha _{s \oplus v} =
 l_{\exp(-2 \pi i s)} r_{\exp(2 \pi i v)}.$$
This action may be viewed as an action of
$H= T \times T$ in the sense of \cite{Rf6}.
For each $\nu$ in the weight lattice $P$ of $\mathfrak g$,
the element $H_{\nu}$ is in $\mathfrak{h}_{\mathbb R}$. We use the notation
$H_{\nu } \oplus H_{\mu}$ to denote $H_{\nu } + H_{\mu}$ as an element
of $ {\mathfrak{h}}_{\mathbb R} \oplus  {\mathfrak{h}}_{\mathbb R} $.
Keep the notation of \cite{Rf6}
for the spectral subspaces of the action $\alpha$ (see 2.22 there).

Let $\check{u}$ be the map on ${\mathfrak h}$*
determined by $u$ via the Killing form
$(\; , \; )$ on $\mathfrak g$.
Let
\bgeqq
& & p= -(H_{\nu_1} \oplus H_{\mu_1}), \; \; \;
q= -(H_{\nu_2} \oplus H_{\mu_2}), \\
& & J=  \frac{h}{4 \pi} (S_u \oplus (-S_u)),
\ndeqq
where $S_u$ is the skew-symmetric operator on $\mathfrak{h}_{\mathbb R}$
defined by
$$S_u (H_\nu) = \sum_{k,l} c_{kl} \nu(H_k) H_l. $$
Then one has
\bgeqq
& &
C^{\lambda_1}_{ \mu_1,i_1; \nu_1, j_1} \circ
C^{\lambda_2}_{  \mu_2,i_2; \nu_2, j_2} \; \\
& & = \exp( \frac{ih}{2}(( {\mu}_1, \check{u} {\mu}_2 )
- ( \nu _1, \check{u} \nu _2 )))
C^{\lambda_1}_{ \mu_1, i_1 ; \nu_1, j_1}
C^{\lambda_2}_{ \mu_2, i_2 ; \nu_2, j_2} \\
& &  = \exp(- 2 \pi i <p, J q>)
C^{\lambda_1}_{ \mu_1, i_1 ; \nu_1, j_1}
C^{\lambda_2}_{ \mu_2, i_2 ; \nu_2, j_2}
\ndeqq
where $\circ$ on the left-hand side is the multiplication
in $A_{K^u_q}$ and the right-hand side is the multiplication in $A_{K_q}$.

On the other hand one has from 2.22 of \cite{Rf6} that
$$
C^{\lambda_1}_{ \mu_1,i_1; \nu_1, j_1} \times_J
C^{\lambda_2}_{  \mu_2,i_2; \nu_2, j_2}
 =
\exp(- 2 \pi i <p, J q>)
C^{\lambda_1}_{ \mu_1,i_1; \nu_1, j_1}
C^{\lambda_2}_{  \mu_2,i_2; \nu_2, j_2}
$$
This means that we have the following result \cite{W4}.
\bgth
\label{rieffel-deformK_q}
The Hopf *-algebras
${\mathcal A}_{K_q^u}$ and $({\mathcal A}_{K_q}, \times_J)$
are isomorphic.
\ndth

That is $K_q^u = (K_q)_J$ in the notation of \thmref{rieffel-deform},
answering Rieffel's question \cite{Rf8} in the affirmative.

Combining \thmref{rieffel-deformK_q}, \thmref{simpleK_q} and
\thmref{simple-invariance}, we obtain the following

\bgth
\label{simpleK^u_q}
Let $K$ be a connected and simply connected simple compact Lie group.
Then for each each $(q, u)$,  $K^u_q$ is an almost classical simple compact quantum group with property $F$.
\ndth

\noindent
{\bf Acknowledgments.}
The research reported here was partially supported by the
National Science Foundation grant DMS-0096136.
During the writing of this work, the author also received support from
the Research Foundation of the University of Georgia in the form of a
Faculty Research Grant in the summer of 2000 and
the Max Planck Institute of Mathematics at Bonn in the form of a visiting
membership for 6 weeks in the summer of 2002. During the summer of 2008,
the author also received a very generous research grant from the mathematics department
of the University of Georgia to enable him to finalize the paper.
The author is grateful to these agencies and institutions
for generous support.


\begin{thebibliography}{99}


\bibitem{Andrus95a} Andruskiewitsch, N., Devoto, J.:
{Extensions of Hopf algebras},
{\em Algebra i Analiz}, {\bf 7}:1, 22--61 (1995).
translation in {\em St. Petersburg Math. J.} {\bf 7}:1,  17--52 (1996).

\bibitem{Banica1} Banica, T.:
{\rm Th\'eorie des repr\'esentations du groupe quantique
compact libre $O(n)$,}
{\em C. R. Acad. Sci. Paris} t. {\bf 322}, Serie I, 241-244 (1996).

\bibitem{Banica2} Banica, T.:
{\rm Le groupe quantique compact libre $U(n)$,}
{\em Commun. Math. Phys.} {\bf 190}, 143-172 (1997).

\bibitem{Banica7} Banica, T.:
Symmetries of a generic coaction,
math.QA/9811060,
{\em Math. Ann.} {\bf 314}, 763-780 (1999).

\bibitem{Banica99c} Banica, T.:
Fusion rules for representations of compact quantum groups,
{\em Exposition Math.} {\bf 17}, 313-337 (1999).

\bibitem{Banica2000a} Banica, T.:
Compact Kac algebras and commuting squares,
{\em J. Funct. Anal.} {\bf 176}, no. 1,
80--99 (2000).

\bibitem{Banica2002}
Banica, T.:
Quantum groups and Fuss-Catalan algebras.
{\em Comm. Math. Phys.} {\bf 226}, no. 1, 221--232 (2002).


\bibitem{Banica2005a}
Banica, T.:
Quantum automorphism groups of homogeneous graphs.
{\em J. Funct. Anal.} {\bf 224}, no. 2, 243--280 (2005).

\bibitem{Banica2005b}
Banica, T.:
Quantum automorphism groups of small metric spaces.
{\em Pacific J. Math.} {\bf 219}, no. 1, 27--51 (2005).

\bibitem{Banica-Bichon2007a}
Banica, T. and  Bichon, J.:
Free product formulae for quantum permutation groups.
{\em J. Inst. Math. Jussieu} {\bf 6}, no. 3, 381--414 (2007).

\bibitem{Banica-Bichon2007b}
Banica, T. and  Bichon, J.:
Quantum automorphism groups of vertex-transitive graphs of order $\leq11$.
{\em J. Algebraic Combin.} {\bf 26}, no. 1, 83--105 (2007).

\bibitem{Banica-Bichon2007prb}
Banica, T. and  Bichon, J.:
Quantum groups acting on 4 points, math/0703118,
to appear in {\em J. Reine Angew. Math.}



\bibitem{Banica-Bichon-Chenevier2007}
Banica, T. and Bichon, J. and Chenevier, G.:
Graphs having no quantum symmetry.
{\em Ann. Inst. Fourier (Grenoble)} {\bf 57}, no. 3, 955--971 (2007).

\bibitem{Banica-Bichon-Collins2007}
Banica, T. and Bichon, J. and Collins, B.:
The hyperoctahedral quantum group.
{\em J. Ramanujan Math. Soc.} {\bf  22},  no. 4, 345--384 (2007).



\bibitem{Bichon1} Bichon, J.:
Quantum automorphism groups of finite graphs.
{\em Proc. Amer. Math. Soc.} {\bf 131}, 665-673 (2003).

\bibitem{Bichon2} Bichon, J.:
Free wreath product by the quantum permutation group.
{\em Alg. Rep. Theory} {\bf 7}, 343-362 (2004).

\bibitem{Bichon3} Bichon, J.:
Algebraic quantum permutation groups.
arXiv:0710.1521. {\em Asian-Eur. J. Math.} {\bf 1}  (2008),  no. 1, 1--13.


\bibitem{Bi-R-V06}
Bichon, J.,  De Rijdt, A., Vaes, S.:
Ergodic coactions with large multiplicity and monoidal equivalence of quantum groups.
{\em Comm. Math. Phys.} {\bf 262} (2006), 703--728.



\bibitem{ChariPressley} Chari, V. and Pressley, A.:
{\em A Guide to Quantum Groups,}
Cambridge University Press, 1994.

\bibitem{Childs} Childs, Lindsay N.:
{\em Taming Wild Extensions: Hopf Algebras and
Local Galois Module Theory,}
American Mathematical Society,  2000.


\bibitem{Connes08mfd} Connes, A.:
On the spectral characterization of manifolds.
arXiv:0810.2088


\bibitem{Dascalescu}
D\u{a}sc\u{a}lescu, S. and   N\u{a}st\u{a}sescu, C. and Raianu, \c{S}.
{\em Hopf Algebras}, Marcel Dekker, Inc., New York, Basel, 2001.

\bibitem{DijKwd1}  Dijkhuizen, M. S. and Koornwinder, T.H.:
{\rm CQG algebras, a direct algebraic approach to compact quantum groups,}
{\em Lett. Math. Phys.} {\bf 32}, 315-330 (1994).

\bibitem{Dr1} Drinfeld, V.G., {Quantum groups},
{\em in} {Proc. of the ICM-1986, Berkeley,} Vol I,
Amer. Math. Soc., Providence, R.I., 1987, pp798--820.



\bibitem{Handel94a} Handelman, David:
Representation ring as invariants for compact groups and limit ratio theorem
for them.
{\em Internat. J. Math.} {\bf 4}, 59-88 (1994).

\bibitem{HR2} Hewitt, E. and Ross, K.:
{\em Abstract Harmonic Analysis}   II,
Springer--Verlag,      1970.

\bibitem{Jimb1} Jimbo, M.,
{A $q$-difference analogue of
$U {\mathfrak g}$ and the Yang-Baxter equations,}
{\em Lett.  Math. Phys.} {\bf 10}, 63-69 (1985).

\bibitem{KlimykSchmudgen} Klimyk, A. U.  and Schm\"{u}dgen, K.:
{\em Quantum Groups and Their Representations},
Springer--Verlag, 1997.

\bibitem{KorogodskiSoib1} Korogodski, Leonid I. and Soibelman, Yan S.,
{\em Algebra of Functions on Compact Quantum Groups}, Part I.
Mathematical Surveys and Monographs, 56. AMS, Providence, RI, 1998

\bibitem{Lev1} Levendorskii, S.:
Twisted algebra of functions on compact quantum group and
their representations,
{\em Algebra i analiz} {\bf 3}:2, 180-198 (1991).
{\em St. Petersburg Math. J.} {\bf 3}:2, 405-423 (1992).

\bibitem{LS1} Levendorskii, S. and Soibelman, Y.:
Algebra of functions on compact quantum groups, Schubert cells,
and quantum tori,
{\em Comm. Math. Phys.} {\bf 139}, 141-170 (1991).

\bibitem{Lus88a} Lusztig, G.:
{\rm Quantum deformations of certain simple modules over enveloping algebras,}
{\em Adv. in Math.} {\bf 70}, 237-249 (1988).


\bibitem{ParshallWang91a} Parshall, B. and Wang, J.:
{\rm Quantum linear groups,}
{\em Memoirs AMS} {\bf 439}, 1991.

\bibitem{Pod6} Podles, P.:
{\rm Symmetries of quantum spaces. Subgroups and
quotient spaces of quantum $SU(2)$ and $SO(3)$ groups,}
{\em Commun. Math. Phys.} {\bf 170}, 1-20 (1995).

\bibitem{PW} Podles, P. and Woronowicz S. L.:
{\rm Quantum deformation of Lorentz group,}
{\em Commun. Math. Phys.} {\bf 130}, 381-431 (1990).

\bibitem{Rf6} Rieffel, M.:
{\rm Deformation quantization for actions of ${\mbox{\bf R}}^d$,}
{\em Memoirs A.M.S.} no. {\bf 506}, 1993.

\bibitem{Rf8} Rieffel, M.:
{\rm Compact quantum groups associated with toral subgroups,}
{\em Contemp. Math.} {\bf 145}, 465-491 (1993).

\bibitem{Rosso87a} Rosso, M.:
{\rm Comparaison des groupes $SU(2)$ quantiques de Drinfeld et Woronowicz,}
{\em C. R. Acad. Sci. Paris} {\bf 304}, 323-326 (1987).

\bibitem{Rosso88a} Rosso, M.:
{\rm Finite dimensional representations of the quantum
analog of the enveloping algebra of a complex semisimple Lie algebra,}
{\em Commun. Math. Phys.} {\bf 117}, 581-593 (1988).

\bibitem{Rosso90a} Rosso, M.:
{\rm Alg\`{e}bres enveloppantes quantifi\'{e}es, groupes quantiques
compacts de matrices et calcul differentiel non-commutatif,}
{\em Duke Math. J.} {\bf 61}, 11-40 (1990).

\bibitem{Schneid93a} Schneider, H.-J.:
Some remarks on exact sequences of quantum groups.
{\em Comm. Algebra} {\bf 21}:9, 3337--3357 (1993).

\bibitem{Schauenburg2000} Schauenburg, P.:
{\rm Faithful flatness over subalgebras: counterexamples.}
{\em Interactions between ring theory and representations of algebras (Murcia)},
331-344, {\em Lecture Notes in Pure and Appl. Math.} {\bf 210}, Dekker, New York, 2000.

\bibitem{Soib1} Soibelman, Y.:
{Algebra of functions on compact quantum group and its representations},
{\em Algebra i analiz} {\bf 2}:1, 190-212 (1990).
{\em Leningrad Math. J.} {\bf 2}:1, 161-178 (1991).

\bibitem{VS}  Soibelman, Y. and Vaksman, L.:
{\rm The algebra of functions on quantum $SU(2)$,}
{\em Funct. Anal. ego Pril.} {\bf  22}:3, 1-14 (1988).
{\em Funct. Anal. Appl.} {\bf 22}:3, 170-181 (1988).

\bibitem{Sweedler} Sweedler, M. E.:
{\em Hopf Algebras,}
Benjamin, New York,  1969.

\bibitem{Takeuchi94a} Takeuchi, Mitsuhiro:
 Quotient spaces for Hopf algebras.
{\em Comm. Algebra} {\bf 22}:7, 2503--2523 (1994).

\bibitem{W5} Van Daele, A. and Wang, S. Z.:
{Universal quantum groups,}
{\em International J. Math.} {\bf 7}, 255-264 (1996).

\bibitem{Wang} Wang, S. Z.:
{\em General Constructions of Compact Quantum Groups,}
Ph.D. Thesis, University of California at Berkeley, March, 1993.

\bibitem{W1} Wang, S. Z.:
{Free products of compact quantum groups,}
{\em Commun. Math. Phys.} {\bf 167}, 671-692 (1995).


\bibitem{W3} Wang, S. Z.:
{\rm Krein duality for compact quantum groups,}
{\em J. Math. Phys.} {\bf 38}:1, 524-534 (1997).

\bibitem{W4} Wang, S. Z.:
{\rm Deformations of compact quantum groups via Rieffel's quantization,}
{\em Commun. Math. Phys.} {\bf 178}:3, 747-764 (1996).

\bibitem{W5'} Wang, S. Z.:
{\rm New classes of compact quantum groups,}
Lecture notes for talks at the University of Amsterdam and
the University of Warsaw, January and March, 1995.

\bibitem{W15} Wang, S. Z.:
{\rm Quantum symmetry groups of finite spaces},
math.OA/9807091,
{\em Commun. Math. Phys.} {\bf 195}:1, 195-211 (1998).

\bibitem{W14} Wang, S. Z.:
{\rm Ergodic actions of universal quantum groups on operator algebras},
{\it Commun. Math. Phys.} {\bf 203}, 481-498 (1999).

\bibitem{W16} Wang, S. Z.:
{\rm Quantum $ax + b$ group as quantum automorphism group of $k[x]$},
{\em Communications in Algebra} {\bf 30}:4, 1807-1815 (2002).

\bibitem{W17} Wang, S. Z.:
{\rm Structure and isomorphic classification of compact quantum groups $A_u(Q)$ and $B_u(Q)$},
{\em J. Operator Theory} {\bf  48},  no. 3, suppl., 573--583  (2002).


\bibitem{normal} Wang, S. Z.:
{\rm Equivalent notions of normal quantum subgroups,
compact quantum groups with properties $F$ and $FD$, and other applications.}
Preprint.


\bibitem{Wor4} Woronowicz, S. L.:
{\rm Twisted $SU(2)$ group. An example of noncommutative
differential calculus,}
{\em Publ. RIMS, Kyoto Univ.} {\bf 23}, 117-181 (1987).

\bibitem{Wor5} Woronowicz, S. L.:
{\rm Compact matrix pseudogroups,}
{\em Commun. Math. Phys.} {\bf 111}, 613-665 (1987).

\bibitem{Wor6} Woronowicz, S. L.:
{\rm Tannaka-Krein duality for compact matrix pseudogroups.
Twisted $SU(N)$ groups,}
{\em Invent. Math.} {\bf 93}, 35-76 (1988).

\end{thebibliography}
\end{document}